\documentclass[12pt]{article}
\usepackage{amsmath,amsfonts,amssymb,amsthm}
\usepackage{color}
\input amssym.def
\begin{document}
\def \Z{\mathbb Z}
\def \C{\mathbb C}
\def \R{\mathbb R}
\def \Q{\mathbb Q}
\def \N{\mathbb N}

\def \L{\mathcal L}
\def \A{{\mathcal{A}}}
\def \g{\mathfrak g}
\def \t{\theta}

\def \l{\ell}
\def \gl{{\frak{gl}}}
\def \<{\langle}
\def \>{\rangle}

\def \span{\mbox{span}}
\def \GL{\mbox{GL}}

\newcommand{\ov}{\overline}

\def \be{\begin{equation}\label}
\def \ee{\end{equation}}
\def \bex{\begin{example}\label}
\def \eex{\end{example}}
\def \bl{\begin{lem}\label}
\def \el{\end{lem}}
\def \bt{\begin{thm}\label}
\def \et{\end{thm}}
\def \bp{\begin{prop}\label}
\def \ep{\end{prop}}
\def \br{\begin{rem}\label}
\def \er{\end{rem}}
\def \bc{\begin{coro}\label}
\def \ec{\end{coro}}
\def \bd{\begin{de}\label}
\def \ed{\end{de}}

\newcommand{\m}{\bf m}
\newcommand{\n}{\bf n}
\newcommand{\nno}{\nonumber}
\newcommand{\nord}{\mbox{\scriptsize ${\circ\atop\circ}$}}
\newtheorem{thm}{Theorem}[section]
\newtheorem{prop}[thm]{Proposition}
\newtheorem{coro}[thm]{Corollary}
\newtheorem{conj}[thm]{Conjecture}
\newtheorem{example}[thm]{Example}
\newtheorem{lem}[thm]{Lemma}
\newtheorem{rem}[thm]{Remark}
\newtheorem{de}[thm]{Definition}
\newtheorem{hy}[thm]{Hypothesis}

\makeatletter \@addtoreset{equation}{section}
\def\theequation{\thesection.\arabic{equation}}
\makeatother \makeatletter

\begin{center}
{\Large \bf  Trigonometric Lie algebras, affine Kac-Moody Lie algebras,
and equivariant quasi modules for vertex algebras}
\end{center}
\begin{center}
{Hongyan Guo$^{a}$
\footnote{Partially supported by China NSF grant (No.11901224)
and the Fundamental Research Funds for
the Central Universities (No.30106220276)},
 Haisheng Li$^{b}$,
Shaobin Tan$^{c}$
\footnote{Partially supported by China NSF grant (No.12131018)}
and Qing Wang$^{c}$
\footnote{Partially supported by
 China NSF grant (Nos.12071385, 12161141001) and
 the Fundamental Research Funds for the Central Universities (No.20720200067)}\\
$\mbox{}^{a}$School of Mathematics and Statistics,
and Hubei Key Laboratory of Mathematical Sciences,
Central China Normal University, Wuhan 430079, China\\
$\mbox{}^{b}$Department of Mathematical Sciences, Rutgers University,
Camden, NJ 08102, USA\\
$\mbox{}^{c}$School of Mathematical Sciences,
 Xiamen University, Xiamen 361005, China\\}
\end{center}

\begin{abstract}
In this paper, we study a family of infinite-dimensional Lie algebras $\widehat{X}_{S}$,
where $X$ stands for the type: $A,B,C,D$, and $S$ is an abelian group,
which generalize the $A,B,C,D$ series of trigonometric Lie algebras.
Among the main results, we identify $\widehat{X}_{S}$ with what are called the covariant algebras of the affine Lie algebra
 $\widehat{\L_{S}}$ with respect to some automorphism groups,
  where $\L_{S}$ is an explicitly defined associative algebra viewed as a Lie algebra.
 We then show that restricted $\widehat{X}_{S}$-modules of level $\l$ naturally correspond to
 equivariant quasi modules for affine vertex algebras related to $\L_{S}$.
 Furthermore, for any finite cyclic group $S$, we
 completely determine the structures of these four families of Lie algebras,
 showing that they are essentially
  affine Kac-Moody Lie algebras of certain types.
\end{abstract}

\section{Introduction}

  Trigonometric Lie algebra (or called Fairlie-Fletcher-Zachos algebra)
is the trigonometric analogue or quantum deformation of the area-preserving
algebra of the torus (cf. \cite{FFZ1}, \cite{FFZ2}, \cite{F}).
In \cite{G-K-L1}, four families of Lie algebras $A_{\bf \hbar}$, $B_{\bf \hbar}$, $C_{\bf \hbar}$, and $D_{\bf \hbar}$,
where ${\bf \hbar}=(\hbar_{1},\ldots,\hbar_{k})\in\R^{k}$,
 were introduced and studied, which generalize the trigonometric sine Lie algebras
 and have extensive applications in the theory of integrable equations and Toda lattices (see \cite{FFZ1}).
 In particular, a connection between these Lie algebras with ${\bf \hbar}$ of certain special values
 and affine Kac-Moody Lie algebras was obtained in \cite{G-K-L1} (cf. \cite{G-K-L2}).
For $\hbar=\hbar_{1}\in\R$ with $k=1$, these Lie algebras were studied in \cite{LTW}
by using vertex algebra theory, and they were realized as what are called the covariant algebras of
the affine Lie algebra $\widehat{\A}$,
where $\A$ is an explicitly defined subalgebra of $\mathfrak{gl}_{\infty}$.
It was proved therein that if $q:=e^{i\hbar}$ is {\em not} a root of unity,
restricted modules for these Lie algebras of level $\l$ are in one-to-one correspondence with
equivariant quasi modules for the affine vertex algebras associated to affine Lie algebra $\widehat{\A}$ of level $\ell$.

In this paper, we continue to study Lie algebras of this type, including the four families with $q$ a root of unity.
More specifically, we study four families of infinite-dimensional Lie algebras
$\widehat{A}_{S}$, $\widehat{B}_{S}$, $\widehat{C}_{S}$, and $\widehat{D}_{S}$,
where $S$ is an abelian group equipped with a linear character $\chi$.
Among the main results, we identify these Lie algebras with the covariant algebras of the affine Lie algebra
 $\widehat{\L_{S}}$ with respect to certain automorphism groups,
  where $\L_{S}$ is an explicitly defined associative algebra viewed as a Lie algebra.
 We then show that restricted modules of level $\l$ for these Lie algebras naturally correspond to
 equivariant quasi modules for affine vertex algebras closely related to $\L_{S}$.
 Furthermore,  assuming that $S$ is a finite cyclic group, we show that
 these Lie algebras are essentially affine Kac-Moody algebras of certain types.

Next, we continue to give a more detailed account of the contents of this paper.
Let $S$ be any abelian group with a linear character $\chi:S\rightarrow\C^{\times}$.
Define $\widehat{A}_{S}$ to be the Lie algebra generated by
$A_{\alpha,m}$ for $\alpha\in S, m\in \Z$
and by ${\bf c}$ (central), subject to relations
 \begin{eqnarray*}
 [A_{\alpha,m},A_{\beta,n}]  =
    (\chi(m\beta-n\alpha)-\chi(n\alpha-m\beta))A_{\alpha+\beta,m+n}
    +m\delta_{\alpha+\beta,0}\delta_{m+n,0}{\bf c}
\end{eqnarray*}
for $\alpha,\beta\in S,\ m,n\in \Z$.
Then follow \cite{G-K-L1} to define Lie algebras $\widehat{B}_{S}$,
$\widehat{C}_{S}$, and $\widehat{D}_{S}$ as the invariant
subalgebras of $\widehat{A}_{S}$ under certain order-$2$ automorphisms.
Among these Lie subalgebras,  Lie algebras $\widehat{B}_{S}$ and $\widehat{C}_{S}$
are shown to be isomorphic.
Note that taking $S=\Z^{k}$ with $\chi({\bf \gamma})=q_{1}^{\gamma_{1}}\cdots q_{k}^{\gamma_{k}}$
 for ${\bf \gamma}=(\gamma_{1},\dots,\gamma_{k})\in \Z^{k}$,
  where $q_{j}=e^{i\hbar_{j}}$ for $1\le j\le k$, we recover the trigonometric Lie algebras.


Note that an associative algebra $\A$ played an important role in \cite{LTW}.
Similarly, we make use of an associative algebra $\L_S$ associated to an abelian group $S$.
Specifically, $\L_{S}$ as a vector space has a basis $\{ L_{\alpha,\beta}\ |\ \alpha,\beta\in S\}$, where
the multiplication is given by
\begin{eqnarray*}
L_{\alpha,\beta}\cdot L_{\mu,\nu}=\delta_{\alpha+\mu,\beta-\nu}L_{\alpha+\mu,\alpha+\nu}
\end{eqnarray*}
for $\alpha,\beta,\mu,\nu\in S$.
Equip $\L_{S}$ with the bilinear form $\<\cdot,\cdot\>$ defined by
\begin{eqnarray*}
\<L_{\alpha,\beta},L_{\mu,\nu}\>=\delta_{\alpha+\mu,0}\delta_{\beta,\nu}\;\;\;\;
\mbox{for}\;\;\alpha,\beta,\mu,\nu\in S,
\end{eqnarray*}
which is non-degenerate, symmetric, and associative. Viewing $\L_{S}$ as a Lie algebra, we
then have an affine Lie algebra $\widehat{\L_{S}}=\L_{S}\otimes\C[t,t^{-1}]+\C{\bf k}$.
As the first result, we show that Lie algebras $\widehat{X}_{S}$ are isomorphic to
 covariant algebras of $\widehat{\L_{S}}$ as we mention next.

Covariant algebras of affine Lie algebras (see \cite{G-K-K}, \cite{Li2}) are closely related to
twisted affine Lie algebras (see \cite{kac}).
For any (possibly infinite-dimensional) Lie algebra $\g$ equipped with a symmetric invariant bilinear
form $\langle\cdot,\cdot\rangle$, we have a (general) affine Lie algebra
$\widehat{\g}=\g\otimes\C[t,t^{-1}]+\C{\bf k}$.
Assume that $\Gamma$ is a group acting on $\g$ as an automorphism group which preserves $\langle\cdot,\cdot\rangle$,
and $\phi:\Gamma\rightarrow\C^{\times}$ is a linear character. Let
$\Gamma$ act on $\widehat{\g}$ as an automorphism group by
$$g{\bf k}={\bf k},\ \ g(a\otimes t^n)=\phi(g)^n(ga\otimes t^n)\ \ \text{ for }g\in \Gamma, a\in \g, n\in \Z.$$
Then we have the $\Gamma$-invariant subalgebra $\widehat{\g}^{\Gamma}$.
In case that $\Gamma$ is a finite cyclic group, $\widehat{\g}^{\Gamma}$ is known as the $\Gamma$-twisted affine Lie algebra.

From the affine Lie algebra $\widehat{\g}$, one can also construct a Lie algebra by covariants as follows.
Let $\Gamma$ and $\phi$ be given as above. In addition, assume that for every $a,b\in \g$,
$$[ga,b]=0 \ \text{ and }\ \<ga,b\>=0\quad \text{ for all but finitely many }g\in \Gamma.$$
Define a new operation $[\cdot,\cdot]_\Gamma$ on $\widehat{\g}$ by
$[{\bf k}, \widehat{\g}]_{\Gamma}=0=[\widehat{\g},{\bf k}]_{\Gamma}$ and
$$[a\otimes t^m,b\otimes t^n]_{\Gamma}
=\sum_{g\in \Gamma} \phi(g)^m([ga,b]\otimes t^{m+n}+m\delta_{m+n,0}\<ga,b\>{\bf k})$$
for $a,b\in \g, m,n\in \Z$. The fact is that the subspace
$$I_\Gamma:={\rm span}\{ (a\otimes t^n)-\phi(g)^n(ga\otimes t^n)\ | \ g\in \Gamma, a\in \g, n\in \Z\}$$
is an ideal of the non-associative algebra $(\widehat{\g},[\cdot,\cdot]_\Gamma)$ and the quotient is a Lie algebra,
which is called the $(\Gamma,\phi)$-covariant algebra
of the affine Lie algebra $\widehat{\g}$ and denoted by $\widehat{\g}[\Gamma]$.
It is known (see \cite{Li2}) that if $\Gamma$ is finite, $\widehat{\g}[\Gamma]$ is canonically isomorphic to
the $\Gamma$-invariant subalgebra $\widehat{\g}^{\Gamma}$.

As the first main result, we show that Lie algebra $\widehat{A}_{S}$
is isomorphic to the $(S,\chi)$-covariant algebra $\widehat{\L_{S}}[S]$
with an explicitly given isomorphism $\psi_A: \widehat{A}_{S}\rightarrow \widehat{\L_{S}}[S]$.
Then we show that $\widehat{B}_{S}$ and $\widehat{D}_{S}$ are isomorphic to
the covariant algebras $\widehat{\L_{S}}[\tilde{S}]$, where $\tilde{S}$ is an automorphism group of
$\widehat{\L_{S}}$, generated by $S$ and a particular order-$2$ automorphism,
with two different linear characters of $\tilde{S}$.
With this realization, by applying a result of \cite{Li2},
we get an equivalence between the category of restricted $\widehat{X}_{S}$-modules of
level $\l\in\C$ and that of equivariant quasi modules for
affine vertex algebras closely related to $\widehat{\L_{S}}$ of level $\ell$.
These generalize the corresponding results of \cite{LTW}.

The second part of this paper is to explore natural connections of $\widehat{X}_{S}$
with affine Kac-Moody Lie algebras. Here, we focus on the case with $S$ a finite cyclic group.
It turns out that it is better to deal with the two cases with $|S|$ even and with $|S|$ odd separately.
When $|S|$ is odd, using the realization of $\widehat{X}_{S}$ as the covariant algebras of the affine
Lie algebra $\widehat{\L_{S}}$, we show that with $|S|=2l+1$,
 $\L_S$ is naturally isomorphic to $\mathfrak{gl}_{2l+1}$ and
Lie algebra $\widehat{A}_{S}$ is isomorphic to
the affine Kac-Moody Lie algebra $\widehat{\mathfrak{gl}}_{2l+1}$.
Then we show that Lie algebra $\widehat{B}_{S}$ (and $\widehat{C}_{S}$)
is isomorphic to the direct sum of the twisted affine Kac-Moody algebra
$\g'(A_{2l}^{(2)})$ with a twisted Heisenberg algebra, divided by a one-dimensional central ideal.
On the other hand, Lie algebra $\widehat{D}_{S}$ is shown to be isomorphic
to the affine Kac-Moody algebra $\g'(B_{l}^{(1)})$.

In the case with $|S|=2l$ (even), we have a Lie algebra homomorphism from $\L_S$ to $\mathfrak{gl}_{2l}$,
but it is neither injective nor surjective. A closely related fact is that Lie algebras $\widehat{X}_{S}$ contain large ideals.
For this case, by making use of the trigonometric basis of classical simple Lie algebras
(cf. \cite{FFZ2}, \cite{Gao}, \cite{G-K-L1}, \cite{G-K-L2}),
we construct a natural surjective Lie algebra homomorphism from $\widehat{A}_{S}$
to the loop algebra $L(\mathfrak{gl}_{l}(\C))$ with the kernel explicitly determined (Theorem \ref{aS}).
Then we construct a surjective Lie algebra homomorphism from
$\widehat{C}_{S}$ to $L(\mathfrak{gl}_{l}(\C))$ for $l$ odd, or to $L(\mathfrak{sp}(l,\C))$ for $l$ even.
We also construct a surjective Lie algebra homomorphism from $\widehat{D}_{S}$ to $L(\mathfrak{o}(l,\C))$.

Note that Lie algebra $\widehat{D}_{S}$ with $S=\Z$ is closely related to the $q$-Virasoro algebra
(see Proposition \ref{qVir}).
A connection between the  $q$-Virasoro algebra with $q$ a root of unity and
affine Kac-Moody Lie algebras of type $B$ or $D$ was obtained in \cite{gltw1, gltw2}.

This paper is organized as follows:
In Section 2, we introduce Lie algebras
$\widehat{X}_{S}$ for any abelian group $S$ equipped with a linear character.
In Section 3, we recall equivariant quasi modules for vertex algebras and covariant algebras of
affine Lie algebras.
In Section 4, we give a realization of $\widehat{X}_{S}$ as the covariant algebras of the affine Lie algebra
$\widehat{\L_{S}}$, and then show the equivalence between the category of restricted
$\widehat{X}_{S}$-modules of level $\l$ and the category of equivariant quasi modules of
 affine vertex algebras related to $\widehat{\L_{S}}$.
 In Section 5, we present a connection of Lie algebras $\widehat{X}_{S}$
 with certain twisted or untwisted affine Kac-Moody Lie algebras.

\section{Lie algebras $\widehat{X}_{S}$}

In this section, to any abelian group $S$ with a linear character $\chi$,
we associate a family of infinite-dimensional Lie algebras
$\widehat{X}_{S}$ where $X$ stands for the type: $A,B,C$, or $D$.

We begin with a generalization of the trigonometric Lie algebra of type $A$ (cf. \cite{G-K-L1}).

\bd{def-AS}
{\em Let $S$ be an additive abelian group with a group homomorphism $\chi: S\rightarrow \C^{\times}$,
i.e., a linear character of $S$.
Define {\em $\widehat{A}_{S}$} to be the Lie algebra
generated by ${\bf c}$ and $A_{\alpha,m}$ for $\alpha\in S,\ m\in \Z$,
subject to the relation that ${\bf c}$ is central and
 \begin{eqnarray}\label{Abracket}
 [A_{\alpha,m},A_{\beta,n}]  =
    (\chi(m\beta-n\alpha)-\chi(n\alpha-m\beta))A_{\alpha+\beta,m+n}
    +m\delta_{\alpha+\beta,0}\delta_{m+n,0}{\bf c}
\end{eqnarray}
for $\alpha,\beta\in S,\ m,n\in \Z$. }
\ed

It is straightforward to show that in fact ${\bf c}$ and $A_{\alpha,n}$ for $\alpha\in S,\ n\in \Z$
form a basis of $\widehat{A}_{S}$ by proving that the vector space $L$ with basis
$\{{\bf c}\}\cup \{ A_{\alpha,n}\ |\ \alpha\in S,\ n\in \Z\}$
and with the bilinear operation given by (\ref{Abracket}) and by $[{\bf c},L]=0=[L,{\bf c}]$ is a Lie algebra.
(Alternatively, this will follow from Theorem \ref{cov-A}.)

The following follows from a straightforward argument:

\bl{aut-AS}
Let $\tau_B,\tau_C,\tau_D$ be the linear endomorphisms of $\widehat{A}_{S}$ uniquely determined by
\begin{eqnarray*}
\begin{aligned}
&\tau_{B}({\bf c})={\bf c},\quad \tau_{B}(A_{\alpha,m})=-(-1)^{m}A_{-\alpha,m},\\
&\tau_{C}({\bf c})={\bf c},\quad\tau_{C}(A_{\alpha,m})=-(-1)^{m}\chi(2\alpha)A_{-\alpha,m},\\
&\tau_{D}({\bf c})={\bf c},\quad\tau_{D}(A_{\alpha,m})=-\chi(2\alpha)A_{-\alpha,m}
\end{aligned}
\end{eqnarray*}
for $\alpha\in S$, $m\in\Z$.
Then $\tau_{B},\tau_{C},\tau_{D}$ are order-$2$ automorphisms of Lie algebra $\widehat{A}_{S}$.
\el

Next, we consider the fixed-point subalgebras of $\widehat{A}_{S}$ under these automorphisms.

\bd{def-ABCD}
{\em For $X=B,C$ or $D$, define $\widehat{X}_{S}$ to be the subalgebra  of the $\tau_{X}$-fixed points in
$\widehat{A}_{S}$, namely the $\tau_X$-invariant subalgebra:
\begin{eqnarray}
\widehat{X}_{S}:=(\widehat{A}_{S})^{\tau_{X}}=\{ u\in \widehat{A}_{S}\ |\ \tau_{X}(u)=u\}.
\end{eqnarray}}
\ed

Among these Lie (sub)algebras, Lie algebras $\widehat{B}_{S}$ and $\widehat{C}_{S}$ are actually isomorphic.
More specifically, using the definition of $\widehat{A}_S$ straightforwardly, we get:

\bl{lemm-B=C}
 Lie algebra $\widehat{A}_S$ admits an automorphism $\sigma$ such that
\begin{eqnarray}
\sigma (A_{\alpha,m})=\chi(\alpha)A_{\alpha,m} \ \ \text{ for }\alpha\in S, m\in \Z.
\end{eqnarray}
Furthermore, we have $\tau_B=\sigma\tau_C \sigma^{-1}$ and $\widehat{B}_{S}\simeq \widehat{C}_{S}$.
\el

The following is another straightforward result:

\bl{lemma-tau}
Lie algebra $\widehat{A}_S$ admits an automorphism $\tau$ such that
\begin{eqnarray}
\tau (A_{\alpha,m})=-A_{-\alpha,m} \ \ \text{ for }\alpha\in S,\; m\in \Z.
\end{eqnarray}
Furthermore, we have $\tau_D=\tau \sigma^2=\sigma^{-1}\tau \sigma$.
\el


Note that for any order-$2$ automorphism $\sigma$ of a Lie algebra $\g$, we have
\begin{eqnarray*}
\g^{\sigma}=\{ u+\sigma(u)\ |\  u\in \g\}.
\end{eqnarray*}
For $\alpha\in S, m\in\Z$, set
\begin{eqnarray}
&&B_{\alpha,m}=A_{\alpha,m}+\tau_{B}(A_{\alpha,m})=A_{\alpha,m}-(-1)^{m}A_{-\alpha,m},\\
&&C_{\alpha,m}=A_{\alpha,m}+\tau_{C}(A_{\alpha,m})=A_{\alpha,m}-(-1)^{m}\chi(2\alpha)A_{-\alpha,m},\\
&&D_{\alpha,m}=A_{\alpha,m}+\tau_{D}(A_{\alpha,m})=A_{\alpha,m}-\chi(2\alpha)A_{-\alpha,m}.
\end{eqnarray}
As $\widehat{B}_{S}\simeq \widehat{C}_{S}$, we shall just study $\widehat{B}_{S}$ or $\widehat{C}_{S}$.

Note that $\widehat{B}_{S}$ is linearly spanned by ${\bf c}$ and $B_{\alpha,m}$ for $\alpha\in S,  m\in\Z$.
It is straightforward to show that the following relations hold for $\alpha,\beta\in S, m,n\in \Z$:
\begin{eqnarray}\label{B-symmetrybracket}
\begin{aligned}
& B_{-\alpha,m}=-(-1)^{m}B_{\alpha,m},\\
 [B_{\alpha,m},B_{\beta,n}]  &=
    (\chi(m\beta -n\alpha)-\chi(n\alpha-m\beta))B_{\alpha+\beta,m+n}
                                                    \\
    &\quad +(-1)^{n}(\chi(m\beta +n\alpha)-\chi(-m\beta-n\alpha))B_{\alpha-\beta,m+n}
                                                               \\
    &\quad +2m(\delta_{\alpha+\beta,0}-(-1)^{m}\delta_{\alpha,\beta})\delta_{m+n,0}{\bf c}.
    \end{aligned}
\end{eqnarray}
Set
\begin{eqnarray}
S^{0}=\{\alpha\in S\ | \ 2\alpha=0 \},
\end{eqnarray}
a subgroup of $S$.

 \bp{generators-relations-BCD}
 Lie algebra $\widehat{B}_{S}$ is isomorphic to the Lie algebra generated by
${\bf c}$ and $B_{\alpha,m}$ for $\alpha\in S,\ m\in \Z$,
subject to relations (\ref{B-symmetrybracket}) and $[{\bf c},\widehat{B}_{S}]=0$.
On the other hand, Lie algebra $\widehat{D}_{S}$ is isomorphic to the Lie algebra generated by set
$\{ D_{\alpha,m}\ |\ \alpha\in S,\ m\in \Z\}\cup \{ {\bf c}\},$
subject to relations $[{\bf c},\widehat{D}_{S}]=0$,
\begin{eqnarray}\label{D-symmetrybracket}
\begin{aligned}
&D_{-\alpha,m}=-\chi(-2\alpha)D_{\alpha,m},\\
 [D_{\alpha,m},D_{\beta,n}] & =
    (\chi(m\beta -n\alpha)-\chi(n\alpha-m\beta))D_{\alpha+\beta,m+n}
                                                    \\
    &\quad +\chi(2\beta)(\chi(m\beta +n\alpha)-\chi(-m\beta-n\alpha))D_{\alpha-\beta,m+n}
                                                               \\
    &\quad +2m(\delta_{\alpha+\beta,0}-\chi(2\alpha)\delta_{\alpha,\beta})\delta_{m+n,0}{\bf c}
    \end{aligned}
\end{eqnarray}
 for $\alpha,\beta\in S,\ m,n\in \Z$.
\ep

\begin{proof} Let $K$ be the Lie algebra with generators
${\bf c}'$ and $B'_{\alpha,m}$ for $\alpha\in S,\ m\in \Z$,
subject to relations (\ref{B-symmetrybracket}) (with $B'_{\alpha,m}$ in place of $B_{\alpha,m}$), and $[{\bf c}',K]=0$.
Then there exists a Lie algebra homomorphism
$\psi: K\rightarrow \widehat{B}_{S}$ such that $\psi({\bf c}')={\bf c}$ and $\psi(B'_{\alpha,m})=B_{\alpha,m}$ for $\alpha\in S,\ m\in \Z$.
It is clear that $\psi$ is onto. It remains to show that $\psi$ is also one-to-one.
Define an equivalence relation $\sim$ on $S$ by $\alpha\sim \beta$ if $\alpha=\beta$ or $\alpha=-\beta$.
Let $T$ be a complete set of equivalence class representatives of $S$.
Note that for $\alpha\in S^0$, we have $[\alpha]=\{\alpha\}\subset S^0$.
Then $S^0\subset T$. Set $T_1=T\backslash S^0$.
Then ${\bf c}$ and $B_{\alpha,m}$ for $(\alpha,m)\in (T_1\times \Z)\cup (S^0\times (2\Z+1))$
form a basis of $\widehat{B}_{S}$. On the other hand, from relations (\ref{B-symmetrybracket}),
${\bf c}'$ and $B_{\alpha,m}'$ for $(\alpha,m)\in (T_1\times \Z)\cup (S^0\times (2\Z+1))$
linearly span $K$.
Then it follows that $\psi$ is one-to-one. For $\widehat{D}_S$, it can be proved similarly.
\end{proof}


\br{nonroot-of-unity}
{\em Let $\hbar_{1},\ldots,\hbar_{k}\in \R$.
Set $q_{r}=e^{i\hbar_{r}}$ for $1\le r\le k$ (where $i=\sqrt{-1}$).
Take $S=\Z^{k}$ with the linear character $\chi$ defined by
$$\chi({\bf n})=q_{1}^{n_{1}}\cdots q_{k}^{n_{k}}
\quad \text{  for }{\bf n}=(n_{1},\ldots,n_{k})\in \Z^{k}.$$
Then $\widehat{X}_{S}$ coincide with the trigonometric Lie algebras
 of type $X$ with $X=A$, $B$, $C$, or $D$, introduced in \cite{G-K-L1}.
For $k=1$, these Lie algebras were studied in \cite{LTW} in the context of
affine vertex algebras and their equivariant quasi modules.}
\er

The Lie algebra $\widehat{D}_{S}$ is closely related to the $q$-Virasoro algebra $D_{S}$ studied in \cite{gltw2}.
Assume $\chi$ is one-to-one.
It was proved therein (Lemma 2.4) that $D_{S}$ is isomorphic to the Lie algebra generated by set
$\{ {\bf k}\}\cup \{ \tilde{D}^{\alpha}(m)\ |\  \alpha\in S, \  m\in \Z\},$
subject to relations
$$[{\bf k},D_{S}]=0,\;\;\;\;
\tilde{D}^{-\alpha}(m)=-\tilde{D}^{\alpha}(m),$$
and
\begin{eqnarray}\label{qVirbracket}
&&[\tilde{D}^{\alpha}(m),\tilde{D}^{\beta}(n)]\nonumber\\
&=&
    (\chi(n\alpha-m\beta)-\chi(m\beta -n\alpha))\tilde{D}^{\alpha+\beta}(m+n)
                                                    \nonumber\\
   &&+(\chi(-n\alpha-m\beta)-\chi(m\beta +n\alpha))\tilde{D}^{\alpha-\beta}(m+n)
                                                               \nonumber\\
   & &+ m\Big(\chi((m+1)(\alpha+\beta))\delta_{2(\alpha+\beta),0}-
             \chi((m+1)(\alpha-\beta))\delta_{2(\alpha-\beta),0}\Big)\delta_{m+n,0}{\bf k}\  \  \  \
\end{eqnarray}
for $\alpha,\beta\in S,\ m,n\in \Z$.
The following is straightforward:

\bp{qVir}
Let $S$ be an additive abelian group such that $S^{0}=0$. Assume $\chi:S\rightarrow \C^{\times}$ is one-to-one.
Then the correspondence
\begin{eqnarray*}
 {\bf c}\mapsto \frac{1}{2}{\bf k}, \ \
D_{\alpha,m}\mapsto -\chi(\alpha)\tilde{D}^{\alpha}(m)\quad \text{ for }\alpha\in S, \ m\in\Z,
\end{eqnarray*}
gives a Lie algebra isomorphism from $\widehat{D}_{S}$ to the $q$-Virasoro algebra $D_{S}$.
\ep

\section{Covariant algebras of affine Lie algebras and quasi modules for affine vertex algebras}

In this section, we recall from \cite{Li1, Li2} the basic notions and results about
$\Gamma$-vertex algebras and their equivariant quasi modules.
In particular, we recall the connection between restricted modules for what were called the covariant algebras of affine Lie algebras
and equivariant quasi modules for affine vertex algebras.

We begin by recalling the notion of quasi module for a vertex algebra, which was introduced in \cite{Li1},
naturally generalizing that of module.

\bd{def-quasi-module}
{\em Let $V$ be a vertex algebra. {\em A quasi $V$-module} is a vector space $W$ equipped with a
linear map
\begin{eqnarray*}
Y_W(\cdot,x): &V&\rightarrow ({\rm End} W)[[x,x^{-1}]]\\
&v&\mapsto Y_W(v,x),
\end{eqnarray*}
satisfying the condition that $Y_W({\bf 1},x)=1$ (the identity operator on $W$),
$$Y_W(v,x)w\in W((x))\quad \text{ for }v\in V,\ w\in W,$$
and for any $u,v\in V$, there exists a nonzero polynomial $p(x_1,x_2)$ such that
\begin{eqnarray}
&&x_0^{-1}\delta\left(\frac{x_1-x_2}{x_0}\right)p(x_1,x_2)Y_W(u,x_1)Y_W(v,x_2)\nonumber\\
&&\quad\quad -x_0^{-1}\delta\left(\frac{x_2-x_1}{-x_0}\right)p(x_1,x_2)Y_W(v,x_2)Y_W(u,x_1)\nonumber\\
&&= x_2^{-1}\delta\left(\frac{x_1-x_0}{x_2}\right)p(x_1,x_2)Y_W(Y(u,x_0)v,x_2).
\end{eqnarray}}
\ed

\bd{def-gamm-va}
{\em Let $\Gamma$ be a group.
A {\em $\Gamma$-vertex algebra} is a vertex algebra $V$ equipped with group homomorphisms
\begin{eqnarray*}
R:\Gamma\rightarrow \GL(V);\; g\mapsto R_{g},
\;\;\mbox{and}\;\;\phi:\Gamma\rightarrow\C^{\times}
\end{eqnarray*}
such that $R_{g}({\bf 1})={\bf 1}$ for $g\in\Gamma$
and
\begin{eqnarray*}
R_{g}Y(v,x)R_{g}^{-1}=Y(R_{g}(v),\phi(g)^{-1}x)\  \  \mbox{for }g\in\Gamma,\  v\in V.
\end{eqnarray*}}
\ed

\br{Z-graded-va-gamma-va}
{\em Let $V$ be a $\Z$-graded vertex algebra in the sense that $V$ is a vertex algebra equipped with a $\Z$-grading
$V=\oplus_{n\in \Z}V_{(n)}$ such that
$$u_pV_{(n)}\subset V_{(n+m-p-1)}\quad \text{ for }u\in V_{(m)},\ m,n,p\in \Z.$$
 Denote by $L(0)$ the degree operator on $V$.
Suppose that $\Gamma$ is a group acting as an automorphism group on $V$ (preserving the grading)
and $\phi:\Gamma\rightarrow\C^{\times}$ is a group homomorphism.
Then $V$ becomes a $\Gamma$-vertex algebra with $R_{g}=\phi(g)^{-L(0)}g$
for $g\in\Gamma$ (see Proposition 6.8 of \cite{Li1}).}
\er

\bd{def-equivariant-module}
{\em Let $V$ be a $\Gamma$-vertex algebra. {\em An equivariant quasi $V$-module}
is a quasi $V$-module $(W,Y_W)$ satisfying the conditions that
\begin{eqnarray}
Y_{W}(R_g(v),x)=Y_W(v,\phi(g)x)\quad \text{ for }g\in \Gamma,\ v\in V
\end{eqnarray}
and that for any $u, v\in V$, there exist $\alpha_1,\dots,\alpha_r\in \phi(\Gamma)\subset \C^{\times}$ such that
\begin{eqnarray}
(x_1-\alpha_1x_2)\cdots (x_1-\alpha_rx_2) [Y_{W}(u,x_1),Y_{W}(v,x_2)]=0.
\end{eqnarray}}
\ed

Next, we discuss covariant algebras of general affine Lie algebras.
Let $\g$ be any (possibly infinite-dimensional) Lie algebra  equipped with a symmetric
invariant bilinear form $\<\cdot,\cdot\>$.  Then we have an (untwisted) affine Lie algebra
\begin{eqnarray*}{\label{eq:aff}}
\widehat{\g}=\g\otimes\C[t,t^{-1}]\oplus\C{\bf k},
\end{eqnarray*}
where ${\bf k}$ is central and
$$[a\otimes t^{m}, b\otimes t^{n}]=[a,b]\otimes t^{m+n}+m\<a,b\>\delta_{m+n,0}{\bf k}$$
for $a,b\in\g$, $m,n\in\Z$.

\br{rem-affine-automorphism}
{\em Assume that $\Gamma$ is a group acting on $\g$ as an automorphism group
which preserves the bilinear form $\<\cdot,\cdot\>$.
 Let $\chi: \Gamma\rightarrow \C^{\times}$ be any linear character.
 Then $\Gamma$ acts on $\widehat{\g}$ as an automorphism group by
 \begin{eqnarray}
 g{\bf k}={\bf k},\ \ g(a\otimes t^n)=\chi(g)^n(ga\otimes t^n)\quad \text{ for }g\in \Gamma, a\in \g,n\in \Z.
 \end{eqnarray}}
 \er

 Recall that for any $\l\in \C$, we have
 the (universal) affine vertex algebra $V_{\widehat{\g}}(\l,0)$ (cf. \cite{FZ}, \cite{LL}),
whose underlying space is the level $\ell$ generalized Verma $\widehat{\g}$-module
$$V_{\widehat{\g}}(\l,0)
=U(\widehat{\g})\otimes_{U(\g\otimes\C[t]+\C{\bf k})}\C_{\l},$$
where $\C_{\l}$ denotes the one-dimensional $(\g\otimes\C[t]+\C{\bf k})$-module
$\C$ with $\g\otimes\C[t]$ acting trivially and with ${\bf k}$ acting as scalar $\l$.
Vertex algebra $V_{\widehat{\g}}(\l,0)$ is naturally a $\Z$-graded vertex algebra
with $V_{\widehat{\g}}(\l,0)_{(0)}=\C {\bf 1}$ and $V_{\widehat{\g}}(\l,0)_{(1)}=\g$.

 Assume that $\Gamma$ is a group acting on $\g$ as an automorphism group
which preserves the bilinear form $\<\cdot,\cdot\>$.
Then $\Gamma$ with the trivial character acts $\widehat{\g}$ as an automorphism group
and then acts on $V_{\widehat{\g}}(\l,0)$
 as an automorphism group preserving the $\Z$-grading.
In view of Remark \ref{Z-graded-va-gamma-va},
given any linear character $\phi:\Gamma\rightarrow\C^{\times}$, we can make $V_{\widehat{\g}}(\l,0)$  a $\Gamma$-vertex algebra.

To introduce covariant algebras of affine Lie algebras, we first recall the following result from \cite{Li2} (Lemma 4.1):

\bl{lemma-covariant}
Let $K$ be a Lie algebra  equipped with a symmetric invariant bilinear form $\<\cdot,\cdot\>$.
Assume that $\Gamma$ is a group acting on $K$ as an automorphism group such that
 $\Gamma$ preserves the bilinear form $\<\cdot,\cdot\>$ and that for $u,v\in K$,
$$[gu,v]=0\;\;\;\mbox{and}\;\;\;\<gu,v\>=0\quad \mbox{for all but finitely many }g\in\Gamma.$$
Define a new bilinear operation $[\cdot,\cdot]_\Gamma$ on $K$ by
\begin{eqnarray}{\label{cov}}
[u,v]_{\Gamma}
=\sum\limits_{g\in\Gamma}[gu,v]
\end{eqnarray}
for $u,v\in K$. Set
$$I_{\Gamma}={\rm span}\{ gu-u\ |\ g\in \Gamma,\ u\in K\}.$$
Then $I_{\Gamma}$ is a two-sided ideal of the non-associative algebra $(K, [\cdot,\cdot]_\Gamma)$
and the quotient algebra $K/I_{\Gamma}$ is a Lie algebra,
which we denote by $K/\Gamma$. Define a bilinear form $\<\cdot,\cdot\>_{\Gamma}$ on $K$ by
\begin{eqnarray}
\<u,v\>_{\Gamma}=\sum_{g\in \Gamma}\<gu,v\>.
\end{eqnarray}
Then $\<\cdot,\cdot\>_{\Gamma}$ reduces to a symmetric invariant bilinear form on $K/\Gamma$.
\el

\br{rem-inv=cova}
{\em Let $K$ be a Lie algebra with a symmetric invariant bilinear form $\<\cdot,\cdot\>$ and
let $\Gamma$ be a finite group acting on $K$ as an automorphism group, preserving $\<\cdot,\cdot\>$.
On the one hand, we have the $\Gamma$-invariant Lie subalgebra $K^{\Gamma}$ and
on the other hand, we have the Lie algebra $K/\Gamma$. Then the linear map $\psi: K\rightarrow K$ defined by
$\psi(u)=\sum_{g\in \Gamma}gu$ gives rise to a Lie algebra isomorphism from $K/\Gamma$ to $K^{\Gamma}$ and we have
$$\<\psi(u),\psi(v)\>=|\Gamma|\<u,v\>_{\Gamma}\ \ \text{for }u,v\in K.$$}
\er

Here, we also have the following result:

\bl{quotient}
Under the setting of Lemma \ref{lemma-covariant}, let $H$ be a normal subgroup of $\Gamma$.
Then $\Gamma$ naturally acts on Lie algebra $K/H$ as an automorphism group
preserving the bilinear form while $H$ acts on $K/H$ trivially.
Furthermore, $K/\Gamma$ is isomorphic to $(K/H)/(\Gamma/H)$.
\el

\begin{proof}  Set $I_H={\rm span}\{ hu-u\ |\ h\in H,u\in K\}$. As $H\lhd \Gamma$,  we have
$$g(hu-u)=(ghg^{-1})gu-gu\in I_H
\quad \text{for }g\in \Gamma, h\in H, u\in K.$$
On the other hand, for $g\in \Gamma, u,v\in K$, we have
$$[gu,gv]_H=\sum_{h\in H}[hgu,gv]=\sum_{h\in H}g[g^{-1}hgu,v]=g[u,v]_H,$$
$$\<gu,gv\>_H=\sum_{h\in H}\<hgu,gv\>=\sum_{h\in H}\<g^{-1}hgu,v\>=\<u,v\>_H.$$
Then $\Gamma$ acts on the Lie algebra $K/H$ as an automorphism group
preserving the bilinear form. It is clear that $H$ acts on $K/H$ trivially.

Since $I_H\subset I_\Gamma$, the identity operator on $K$ reduces to a surjective linear map
$\theta: K/H\rightarrow K/\Gamma$.
Write $\Gamma=\cup_{j\in J}Hg_j$ disjoint. For $u\in K,\ h\in H,\ j\in J$, we have
$$hg_ju-u=h(g_ju)-g_ju+g_ju-u.$$
It follows that $\theta$ reduces to a linear isomorphism  $\bar{\theta}: (K/I_H)/(\Gamma/H)\rightarrow K/I_{\Gamma}$.
On the other hand,  for $u,v\in K$, we have
 $$[u,v]_{\Gamma}=\sum_{j\in J,h\in H}[hg_ju,v]=\sum_{j\in J}[g_ju,v]_H=([u,v]_H)_{\Gamma/H}.$$
 Then $\bar{\theta}$ is a Lie algebra isomorphism.
\end{proof}

Now, consider the case with $K$ taken to be the affine Lie algebra $\widehat{\g}$.
Assume that $\Gamma$ is a group acting on $\g$ as an automorphism group and $\Gamma$
 preserves the bilinear form $\<\cdot,\cdot\>$.
In addition, let $\phi:\Gamma\rightarrow\C^{\times}$ be a linear character.
Set
\begin{eqnarray}
I_{\Gamma}={\rm span}\{\phi(g)^{m}(ga\otimes t^{m})-(a\otimes t^{m})\ | \ g\in\Gamma,\  a\in\g,\  m\in\Z\},
\end{eqnarray}
a subspace of $\widehat{\g}$.
The following result can be found in  \cite{Li2} [Proposition 4.4] (cf. \cite{G-K-K}):

\bp{cov}
Let $\g$, $\Gamma$ and $\phi$ be given as above.
Assume that for $a,b\in\g$,
$$[ga,b]=0\;\;\;\mbox{and}\;\;\;\<ga,b\>=0\quad \mbox{for all but finitely many}\;g\in\Gamma.$$
Define a new bilinear operation $[\cdot,\cdot]_\Gamma$ on $\widehat{\g}$ by
$[{\bf k},\widehat{\g}]_{\Gamma}=0=[\widehat{\g},{\bf k}]_{\Gamma}$ and
\begin{eqnarray}{\label{eq-cov}}
[a\otimes t^{m},b\otimes t^{n}]_{\Gamma}
=\sum\limits_{g\in\Gamma}\phi(g)^{m}\left([ga,b]\otimes t^{m+n}+m\<ga,b\>\delta_{m+n,0}{\bf k}\right)
\end{eqnarray}
for $a,b\in\g, \  m,n\in\Z$.
Then $I_{\Gamma}$ is a two-sided ideal of the non-associative algebra
 $(\widehat{\g}, [\cdot,\cdot]_\Gamma)$
and the quotient algebra $\widehat{\g}/I_{\Gamma}$ is a Lie algebra.
\ep

Denote the Lie algebra $\widehat{\g}/I_{\Gamma}$ obtained in Proposition \ref{cov} by $\widehat{\g}[\Gamma]$,
which is called the {\em $(\Gamma,\phi)$-covariant algebra of the affine Lie algebra $\widehat{\g}$},
or simply the {\em $\Gamma$-covariant algebra of $\widehat{\g}$} whenever it is clear from the context.
 The image of ${\bf k}$ in $\widehat{\g}[\Gamma]$ is nonzero, so we identify ${\bf k}$ with its image.
 Let $W$ be a $\widehat{\g}[\Gamma]$-module.
For $a\in \g$, set
$$a_W(x)=\sum_{n\in \Z}a(n)_{\Gamma}\; x^{-n-1}\in \mbox{End} (W)[[x,x^{-1}]],$$
where $a(n)_{\Gamma}$ denotes the image of $a\otimes t^n$ in $\widehat{\g}[\Gamma]$.

\br{rcvalgebra}
{\em It is straightforward to show that $\widehat{\g}[\Gamma]$ is
isomorphic to the Lie algebra generated by  $\widehat{\g}$ as a vector space, subject to the following relations:
$[{\bf k},  \widehat{\g}]=0$,
\begin{eqnarray*}
\begin{aligned}
& ga\otimes t^n=\phi(g)^{-n}(a\otimes t^n),\\
[a\otimes t^m,b\otimes t^n] &=\sum_{g\in \Gamma}\phi(g)^m\left([ga,b]\otimes t^{m+n}+m\<ga,b\>\delta_{m+n,0}{\bf k}\right)
\end{aligned}
\end{eqnarray*}
for $g\in \Gamma,\ a,b\in \g,\ m,n\in \Z$.}
\er

The following is an immediate consequence:

\bc{cvla-rep}
Let $W$ be any vector space. Then a linear map
$\rho:  \widehat{\g}\rightarrow {\rm End}(W)$ reduces to
a representation of $\widehat{\g}[\Gamma]$ on $W$ if and only if $[\rho({\bf k}),\rho(\widehat{\g})]=0$,
\begin{eqnarray*}
\rho(ga\otimes t^{n})=\phi(g)^{-n}\rho(a\otimes t^{n})
\end{eqnarray*}
for $g\in \Gamma,\ a\in \g,\ n\in \Z$, and
\begin{eqnarray*}
[\rho(a\otimes t^{m}),\rho(b\otimes t^{n})]
=\sum_{g\in \Gamma}\phi(g)^m\left(\rho([ga,b]\otimes t^{m+n})+m\<ga,b\>\delta_{m+n,0}\rho({\bf k})\right)
\end{eqnarray*}
for $a,b\in \g,\ m, n\in \Z$.
\ec

 \br{rem-g/H}
 {\em Under the setting of Proposition \ref{cov}, let $H$ be a normal subgroup of $\Gamma$.
 In view of Lemma \ref{quotient}, $\Gamma$ acts on Lie algebra $\g/H$
 as an automorphism group (with $H$ acting trivially), preserving the bilinear form.
Then for any linear character $\phi: \Gamma\rightarrow \C^{\times}$,
$V_{\widehat{\g/H}}(\ell,0)$ becomes a $\Gamma$-vertex algebra. }
\er

The following is a result of \cite{Li2} (Theorem 4.9):

\bt{thm-4.9}
Let $\g$ be a Lie algebra  equipped with a symmetric invariant bilinear form $\<\cdot,\cdot\>$.
Assume that $\Gamma$ is a subgroup of ${\rm Aut}(\g,\<\cdot,\cdot\>)$ such that for any $a,b\in \g$,
$$[ga,b]=0\;\;\;\mbox{and}\;\;\;\<ga,b\>=0\quad \mbox{for all but finitely many}\;g\in\Gamma.$$
Let $\phi: \Gamma \rightarrow \C^{\times}$ be a linear character and set $H=\ker \phi\lhd\Gamma$.
Assume $\ell\in \C$.
Then any restricted module $W$ of level $\ell$ for the Lie algebra $\widehat{\g}[\Gamma]$
is a $(\Gamma,\phi)$-equivariant quasi $V_{\widehat{\g/H}}(\ell,0)$-module
 with $Y_W(\bar{a},x)=a_W(x)$ for $a\in \g$, where $\bar{a}$ denotes the image of $a$ in $\g/H$.
 On the other hand, any $(\Gamma,\phi)$-equivariant quasi $V_{\widehat{\g/H}}(\ell,0)$-module $(W,Y_W)$
is a restricted module of level $\ell$ for the Lie algebra $\widehat{\g}[\Gamma]$ with $a_W(x)=Y_W(\bar{a},x)$ for $a\in \g$.
\et

We formulate the following straightforward result:

\bl{lem-affine-covar-invar}
Let $\g$ be a Lie algebra  equipped with a symmetric invariant bilinear form $\<\cdot,\cdot\>$.
Assume that $\Gamma$ is a finite group acting on $\g$ as an automorphism group, preserving
the bilinear form $\<\cdot,\cdot\>$, and assume that  $\phi: \Gamma\rightarrow \C^{\times}$ is a linear character.
Then the linear map $\Psi: \widehat{\g}\rightarrow \widehat{\g}$ defined by
 \begin{eqnarray}
 \Psi({\bf k})=|\Gamma|{\bf k}\ \text{ and }\ \Psi(a\otimes t^n)=\sum_{g\in \Gamma}\phi(g)^n(ga\otimes t^n)
 \ \ \text{ for }a\in \g,\ n\in \Z,
 \end{eqnarray}
 reduces to a Lie algebra isomorphism from $\widehat{\g}[\Gamma]$ to $\widehat{\g}^{\Gamma}$.
\el

As an immediate consequence of Lemma \ref{quotient} we have:

\bc{tensor-covariant}
 Let $\g, \Gamma, \phi$ be given as in Theorem \ref{thm-4.9}. Suppose $\Gamma=\Gamma_1\times \Gamma_2$. Then
 \begin{eqnarray}
 \widehat{\g}[\Gamma]=(\widehat{\g}[\Gamma_1])/\Gamma_2.
 \end{eqnarray}
  On the other hand, $\Gamma$ acts on $\widehat{\g}[\Gamma_1]$ as an automorphism group with $\Gamma_1$ acting trivially.
 \ec



\section{Lie algebras $\widehat{X}_{S}$ and covariant affine Lie algebras}

In this section, given any abelian group $S$, first we introduce an associative algebra $\L_S$
on which $S$ acts naturally as an automorphism group and we equip $\L_S$
with a non-degenerate symmetric associative bilinear form which is preserved by the action of $S$.
Then we realize Lie algebras $\widehat{X}_{S}$ as the covariant algebras of affine Lie algebra
$\widehat{\L_S}$.

The following serves as the definition of the associative algebra $\L_S$:

\bl{assoc-algebra}
Let $S$ be any additive abelian group. Denote by $\L_S$ the vector space with a designated basis
$\{L_{\alpha,\beta}\ | \ \alpha,\beta\in S\}.$
Define a bilinear operation $*$ on $\L_S$ by
\begin{eqnarray}
L_{\alpha,\beta}* L_{\mu,\nu}=\delta_{\alpha+\mu,\beta-\nu}L_{\alpha+\mu,\alpha+\nu}
\;\;\;\;\mbox{for}\;\;\alpha,\beta,\mu,\nu\in S.
\end{eqnarray}
Then $\L_S$ is an associative algebra.
Furthermore, the bilinear form $\<\cdot,\cdot\>$, defined by
\begin{eqnarray}\label{eq:LS-bil}
\<L_{\alpha,\beta},L_{\mu,\nu}\>=\delta_{\alpha+\mu,0}\delta_{\beta,\nu}\quad \text{ for }\alpha,\beta,\mu,\nu\in S,
\end{eqnarray}
 is symmetric, non-degenerate, and associative in the sense that
$\<ab,c\>=\<a,bc\>$ for $a,b,c\in\L_{S}$.
\el

\begin{proof} 
The defined operation $*$ is associative because for $\alpha,\beta,\mu,\nu,p,q\in S$,
\begin{eqnarray*}
L_{\alpha,\beta}* (L_{\mu,\nu}*L_{p,q})=\delta_{\mu+p,\nu-q}L_{\alpha,\beta}*L_{\mu+p,\mu+q}
=\delta_{\mu+p,\nu-q}\delta_{\alpha+\mu+p,\beta-\mu-q}L_{\alpha+\mu+p,\alpha+\mu+q},
\end{eqnarray*}
\begin{eqnarray*}
(L_{\alpha,\beta}* L_{\mu,\nu})*L_{p,q}=\delta_{\alpha+\mu,\beta-\nu}L_{\alpha+\mu,\alpha+\nu}*L_{p,q}
=\delta_{\alpha+\mu,\beta-\nu}\delta_{\alpha+\mu+p,\alpha+\nu-q}L_{\alpha+\mu+p,\alpha+\mu+q},
\end{eqnarray*}
and
$$\delta_{\mu+p,\nu-q}\delta_{\alpha+\mu+p,\beta-\mu-q}=\delta_{\mu+p,\nu-q}\delta_{\alpha+\mu+p,\beta-\nu+p}
=\delta_{\alpha+\mu+p,\alpha+\nu-q}\delta_{\alpha+\mu,\beta-\nu}.$$
It is clear that the defined bilinear form is symmetric and non-degenerate.
Define a linear functional $f$ on $\L_S$ by $f(L_{\alpha,\beta})=\delta_{\alpha,0}$ for $\alpha,\beta\in S$. Then
$$\<L_{\alpha,\beta},L_{\mu,\nu}\>=\delta_{\alpha+\mu,0}\delta_{\beta,\nu}=f(L_{\alpha,\beta}*L_{\mu,\nu}).$$
It follows immediately that the bilinear form is associative.
\end{proof}

We also have the following straightforward results:

\bl{locally-finite}
Group $S$ acts on $\L_{S}$ as an automorphism group with $\gamma\in S$ acting as $\sigma_{\gamma}$ defined by
\begin{eqnarray}
\sigma_{\gamma}(L_{\alpha,\beta})=L_{\alpha,\beta+\gamma}\;\;\;\;\;\mbox{for}\;\;\alpha,\beta,\gamma\in S.
\end{eqnarray}
Furthermore, $S$ preserves the bilinear form $\<\cdot,\cdot\>$.
\el

\bl{tau-D}
Let $\chi: S\rightarrow \C^{\times}$ be a linear character. Define a linear endomorphism $\sigma_{\chi}$ of $\L_S$ by
\begin{eqnarray}
\sigma_{\chi}(L_{\alpha,\beta})=\chi(2\alpha)L_{\alpha,\beta}\ \ \text{ for }\alpha,\beta\in S.
\end{eqnarray}
Then  $\sigma_{\chi}$ is an algebra automorphism of $\L_S$ and preserves the bilinear form.
\el

Viewing $\L_S$ as a Lie algebra with the commutator as its Lie bracket, we have:

\bl{}
The bilinear form $\<\cdot,\cdot\>$ on Lie algebra $\L_{S}$ is symmetric, non-degenerate, and invariant.
Furthermore, for any $\alpha,\beta,\mu,\nu\in S$,
\begin{eqnarray}
[L_{\alpha,\beta+\gamma},L_{\mu,\nu}]=0
\;\;\mbox{and}\;\;\<L_{\alpha,\beta+\gamma},L_{\mu,\nu}\>=0
\end{eqnarray}
for all but finitely many $\gamma\in S$.
\el

Associated to the pair $(\L_S,\<\cdot,\cdot\>)$, we have an affine Lie algebra
\begin{eqnarray*}
\widehat{\L_{S}}=\L_{S}\otimes\C[t,t^{-1}]\oplus\C {\bf k}.
\end{eqnarray*}
Let $\chi$ be a linear character of $S$ as in the definition of Lie algebra $\widehat{A}_S$.
Then we have the covariant Lie algebra $\widehat{\L_{S}}[S]$.
The following is the key result of this section:

\bt{cov-A}
Define a linear map $\psi_A:  \widehat{A}_{S}\rightarrow \widehat{\L_{S}}[S]$ by
\begin{eqnarray*}
\psi_A({\bf c})={\bf k}\ \text{ and }\
 \psi_A(A_{\alpha,m})=L_{\alpha,0}(m)_{S}\quad \text{ for }\alpha\in S, \  m\in\Z,
\end{eqnarray*}
where for any $a\in\L_{S}$, $m\in\Z$, $a(m)_{S}$ denotes the image of $a\otimes t^{m}$ in $\widehat{\L_{S}}[S]$.
Then $\psi_A$ is a Lie algebra isomorphism from $\widehat{A}_{S}$ to the $(S,\chi)$-covariant algebra
 $\widehat{\L_{S}}[S]$ of the affine Lie algebra $\widehat{\L_{S}}$.
\et

\begin{proof}
Recall that  $\{A_{\alpha,m} \ | \ \alpha\in S, \; m\in\Z\}\cup\{ {\bf c}\}$
 is a basis
of $\widehat{A}_{S}$.
On the other hand, from the construction (Proposition \ref{cov}),
$\widehat{\L_{S}}[S]$ as a vector space is the quotient space of $\widehat{\L_{S}}$,
modulo the subspace linearly spanned by vectors
$$\sigma_{\gamma}(a)\otimes t^{m}-\chi(-m\gamma)(a\otimes t^{m})$$
for $a\in\L_{S}$, $\gamma\in S$, $m\in\Z$.
We see that ${\bf k}$ and $L_{\alpha,0}(m)_{S}$ ($\alpha\in S,m\in\Z$) form a basis
of $\widehat{\L_{S}}[S]$.
Let $\alpha,\beta\in S$, $m,n\in\Z$.  For $\gamma\in S$, we have
$$[\sigma_{\gamma}(L_{\alpha,0}), L_{\beta,0}]
=\delta_{\gamma,\alpha+\beta}L_{\alpha+\beta,\alpha}-\delta_{\gamma,-\alpha-\beta}L_{\alpha+\beta,-\alpha},$$
$$\<\sigma_{\gamma}(L_{\alpha,0}), L_{\beta,0}\>
=\delta_{\alpha+\beta,0}\delta_{\gamma,0}.$$
Then we get
\begin{eqnarray}
&&{}[L_{\alpha,0}(m)_{S},L_{\beta,0}(n)_{S}]    \nonumber\\
&&{}=\chi(m(\alpha+\beta))L_{\alpha+\beta,\alpha}(m+n)_{S}
     - \chi(-m(\alpha+\beta))L_{\alpha+\beta,-\alpha}(m+n)_{S}
      +m\delta_{\alpha+\beta,0}\delta_{m+n,0}{\bf k}      \nonumber\\
&&{}=(\chi(m\beta-n\alpha)-\chi(n\alpha-m\beta))
    L_{\alpha+\beta,0}(m+n)_{S}
      +m\delta_{\alpha+\beta,0}\delta_{m+n,0}{\bf k},
\end{eqnarray}
noticing that
$$L_{\alpha+\beta,\alpha}(m+n)_{S}
=\sigma_{\alpha}(L_{\alpha+\beta,0})(m+n)_{S}
=\chi(-(m+n)\alpha)L_{\alpha+\beta,0}(m+n)_{S},$$
$$L_{\alpha+\beta,-\alpha}(m+n)_{S}
=\sigma_{-\alpha}(L_{\alpha+\beta,0})(m+n)_{S}
=\chi((m+n)\alpha)L_{\alpha+\beta,0}(m+n)_{S}.$$
Now it follows from (\ref{Abracket}) that $\psi_{A}$ is a Lie algebra isomorphism.
\end{proof}

Next, we consider the Lie algebra automorphism of $\L_{S}$ corresponding to the automorphism $\tau_B$ of
$\widehat{A}_S$, which defines Lie algebra $\widehat{B}_{S}$.
Note that the linear endomorphism $T$ of $\L_S$
with $T(L_{\alpha,\beta})=L_{-\alpha,\beta}$ for $\alpha,\beta\in S$ is an algebra anti-automorphism of $\L_S$.
It follows that $-T$ is a Lie algebra automorphism of $\L_S$.

\bd{}
{\em Let $\tau$ be the order-2 automorphism of Lie algebra $\L_{S}$ defined by
\begin{eqnarray}
\tau(L_{\alpha,\beta})=-L_{-\alpha,\beta}\;\;\;\;\mbox{for}\;\;\alpha,\beta\in S.
\end{eqnarray}}
\ed

Set
\begin{eqnarray}
\tilde{S}=\<\tau, S\>\subset {\rm Aut} (\mathcal{L}_{S}).
\end{eqnarray}
It is straightforward to see that $\tilde{S}=\<\tau\>\times S$ and $\tilde{S}$ preserves the bilinear form on $\mathcal{L}_{S}$.
Then given any group homomorphism $\phi: \tilde{S}\rightarrow \C^{\times}$,
we have an action of $(\tilde{S},\phi)$ on $\widehat{\L_S}$ as an automorphism group.
Furthermore, we have a covariant algebra $\widehat{\L_S}[\tilde{S}]$.

Define an automorphism $\hat{\tau}$ of the affine Lie algebra $\widehat{\L_S}$ by $\hat{\tau}({\bf k})={\bf k}$ and
\begin{eqnarray}
\hat{\tau}(L_{\alpha,\beta}\otimes t^n)=(-1)^n(\tau(L_{\alpha,\beta})\otimes t^n)=-(-1)^n(L_{-\alpha,\beta}\otimes t^n)
\end{eqnarray}
for $\alpha,\beta\in S, n\in \Z$. It can be readily seen that $\hat{\tau}$ commutes with $S$ on $\widehat{\L_S}$.
Set
\begin{eqnarray}
\tilde{S}_B=\<S,\hat{\tau}\>=S\times \<\hat{\tau}\>\subset {\rm Aut}(\widehat{\L_S}).
\end{eqnarray}
By Lemma \ref{quotient}, we can and we should view $\hat{\tau}$ as an automorphism of $\widehat{\L_S}[S]$.

Recall that $\widehat{B}_S=(\widehat{A}_S)^{\tau_B}$, where
 $$\tau_B(A_{\alpha,n})=-(-1)^nA_{-\alpha,n}\ \  \text{ for }\alpha\in S,n\in \Z,$$
and recall from Theorem \ref{cov-A} the Lie algebra isomorphism $\psi_A:  \widehat{A}_{S}\rightarrow \widehat{\L_{S}}[S]$.
It is straightforward to show that
\begin{eqnarray}
\psi_A\circ \tau_B=\hat{\tau}\circ \psi_A.
\end{eqnarray}
From this we get
\begin{eqnarray}
(\widehat{A}_S)/\<\tau_B\>\overset{\psi_A}{\simeq} (\widehat{\L_S}[S])/\<\hat{\tau}\>.
\end{eqnarray}
 On the other hand, in view of Remark \ref{rem-inv=cova}, there is a Lie algebra isomorphism
 $$\eta: (\widehat{A}_S)^{\tau_B}\rightarrow  \widehat{A}_S/\<\tau_B\>$$
 such that $\eta(B_{\alpha,n})=\overline{B_{\alpha,n}}$ for $\alpha\in S,n\in \Z$, where
 $B_{\alpha,n}=A_{\alpha,n}+\tau_B(A_{\alpha,n})$.

Define a linear character $\tilde{\chi}_B:\tilde{S}_B\rightarrow\C^{\times}$ by
\begin{eqnarray}
\tilde{\chi}_B(\hat{\tau})=-1\; \ \mbox{and}\; \
\tilde{\chi}_B(\gamma)=\chi(\gamma)\; \ \mbox{for} \ \gamma\in S.
\end{eqnarray}
By Lemma \ref{quotient}, we have
$\widehat{\L_S}[\tilde{S}_{B}]\simeq (\widehat{\L_S}[S])/\<\hat{\tau}\>$.
Consequently, we obtain
\begin{eqnarray}
\widehat{B}_S=(\widehat{A}_S)^{\tau_B}\overset{\eta}{\simeq} \widehat{A}_S/\<\tau_B\>\overset{\psi_A}{\simeq} \widehat{\L_S}[S]/\<\hat{\tau}\>
\simeq \widehat{\L_S}[\tilde{S}_B].
\end{eqnarray}
More explicitly,  we have proved (cf. Proposition 2.10 of \cite{LTW}):

\bp{cov-B}
There exists a Lie algebra isomorphism $\Theta$ from $\widehat{B}_{S}$ to the $(\tilde{S}_B,\tilde{\chi}_B)$-covariant algebra
 $\widehat{\L_{S}}[\tilde{S}_B]$, which is uniquely determined by
$$\Theta({\bf c})=\frac{1}{2}{\bf k},\  \
 \Theta(B_{\alpha,m})=L_{\alpha,0}(m)_{\tilde{S}_B}\quad \text{  for }\alpha\in S, \ m\in\Z.$$
 \ep

Next, we consider Lie algebra $\widehat{D}_{S}$.  Recall that $\widehat{D}_{S}=(\widehat{A}_S)^{\tau_D}$, where
$$\tau_D(A_{\alpha,n})=-\chi(2\alpha)A_{-\alpha,n}\quad \text{ for }\alpha\in S,n\in \Z.  $$
Recall from Lemma \ref{tau-D} the automorphism $\sigma_{\chi}$ of associative algebra $\L_S$.
Set $\tau_{\chi}=\tau \sigma_{\chi}$, which is a Lie algebra automorphism of $\L_S$ and preserves the bilinear form. More explicitly, we have
\begin{eqnarray}
\tau_{\chi}(L_{\alpha,\beta})=-\chi(2\alpha)L_{-\alpha,\beta}\quad \text{ for }\alpha,\beta\in S.
\end{eqnarray}
Define an automorphism $\hat{\tau}_{\chi}$ of the affine Lie algebra $\widehat{\L_S}$ by $\hat{\tau}_{\chi}({\bf k})={\bf k}$ and
\begin{eqnarray}
\hat{\tau}_{\chi}(L_{\alpha,\beta}\otimes t^n)=\tau_{\chi}(L_{\alpha,\beta})\otimes t^n
=-\chi(2\alpha)(L_{-\alpha,\beta}\otimes t^n)
\end{eqnarray}
for $\alpha,\beta\in S, n\in \Z$. It is straightforward to show that
\begin{eqnarray}
\psi_A\circ \tau_D=\hat{\tau}_{\chi}\circ \psi_A.
\end{eqnarray}
Set
\begin{eqnarray}
\tilde{S}_D=\<S,\hat{\tau}_{\chi}\>=S\times \<\hat{\tau}_{\chi}\>\subset {\rm Aut}(\widehat{\L_S}),
\end{eqnarray}
and define a linear character $\tilde{\chi}_{D}: \tilde{S}_D\rightarrow\C^{\times}$ by
\begin{eqnarray}
\tilde{\chi}_{D}(\hat{\tau}_{\chi})=1\; \ \mbox{and}\; \
\tilde{\chi}_{D}(\gamma)=\chi(\gamma)\; \ \mbox{for} \ \gamma\in S.
\end{eqnarray}
Similarly, we have
$$\widehat{D}_S=(\widehat{A}_S)^{\tau_D}\simeq \widehat{A}_S/\<\tau_D\>\overset{\psi_A}{\simeq} \widehat{\L_S}[S]/\<\hat{\tau}_{\chi}\>
\simeq \widehat{\L_S}[\tilde{S}_D].$$
Thus we have (cf. Proposition 2.12 of \cite{LTW}):

\bp{cov-D}
The correspondence
$${\bf c}\mapsto \displaystyle\frac{1}{2}{\bf k},\ \
 D_{\alpha,m}\mapsto L_{\alpha,0}(m)_{\tilde{S}_D}\quad \text{ for }\alpha\in S, \ m\in\Z,$$
 gives a Lie algebra isomorphism from $\widehat{D}_{S}$  to the $(\tilde{S}_D,\tilde{\chi}_{D})$-covariant algebra
 $\widehat{\L_{S}}[\tilde{S}_D]$.
\ep

Next, we give a slightly different realization of $\widehat{D}_{S}$.
Recall the automorphism $\tau$ of Lie algebra $\L_{S}$.
With $\tau_D=\sigma^{-1}\tau \sigma$ on $\widehat{A}_S$ (Lemma \ref{lemma-tau}), $\sigma$ is a Lie algebra isomorphism from
$(\widehat{A}_S)^{\tau_D}$ to $(\widehat{A}_S)^{\tau}$.
Then using Lie algebra isomorphism $\psi_A: \widehat{A}_S\rightarrow \widehat{\L_S}[S]$ we get
\begin{eqnarray}\label{DS-isomorphism}
\widehat{D}_S= (\widehat{A}_S)^{\tau_D}\overset{\sigma}{\simeq} (\widehat{A}_S)^{\tau}
\overset{\psi_A}{\simeq}\widehat{\L_S}[S]^{\tau}
\simeq \widehat{\L_S^{\tau}}[S].
\end{eqnarray}
For $\alpha,\beta\in S$, set
 \begin{eqnarray}\label{eq:2.27}
 L^{\tau}_{\alpha,\beta}=L_{\alpha,\beta}+\tau(L_{\alpha,\beta})
 =L_{\alpha,\beta}-L_{-\alpha,\beta}\in \L_S^{\tau}.
 \end{eqnarray}
 For $\alpha\in S, m\in\Z$, we have
 \begin{eqnarray*}
&&\psi_A\sigma(D_{\alpha,m})=\psi_A\sigma(A_{\alpha,m}-\chi(2\alpha)A_{-\alpha,m})
=\chi(\alpha)\psi_A(A_{\alpha,m}-A_{-\alpha,m})\\
&=&\chi(\alpha)(L_{\alpha,0}(m)_S-L_{-\alpha,0}(m)_S)=\chi(\alpha)L_{\alpha,0}^{\tau}(m)_S.
\end{eqnarray*}
 Writing the isomorphism relation (\ref{DS-isomorphism}) explicitly we get (cf. Proposition 2.13 of \cite{LTW}):

\bp{cov-D2}
There is a Lie algebra isomorphism $\Psi$ from $\widehat{D}_{S}$ to  the $(S,\chi)$-covariant algebra
 $\widehat{\L^{\tau}_{S}}[S]$ such that $\Psi({\bf c})={\bf k}$
 and $$\Psi(D_{\alpha,m})=\chi(\alpha)L_{\alpha,0}^{\tau}(m)_{S}\ \ \text{ for }\alpha\in S, m\in\Z.$$
 \ep

Next, we give a precise description of the relationship between Lie algebras $\widehat{X}_{S}$
and affine vertex algebras in terms of equivariant quasi modules.

\begin{de}
{\em
An $\widehat{X}_{S}$-module $W$ is said to be {\em restricted} if for any $\alpha\in S$
and $w\in W$, $X_{\alpha,m}w=0$ for all sufficiently large integers $m$.
We say an $\widehat{X}_{S}$-module $W$ is of {\em level} $\ell\in\C$ if the central element ${\bf c}$
acts as scalar $\l$.}
\end{de}

For $X\in \{A,B,C,D\},\ \alpha\in S$, form a generating function
$$X_{\alpha}(x)=\sum_{m\in\Z}X_{\alpha,m}x^{-m-1}.$$


Let $H$ be a subgroup of $\tilde{S}$.
Then $H$ acts on Lie algebra $\L_S$ as an automorphism group and preserves the bilinear form.
By Lemma \ref{lemma-covariant} we have a Lie algebra $\L_S/H$ with a symmetric invariant bilinear form
given by
$$\<\bar{u},\bar{v}\>=\sum_{h\in H}\<hu,v\>\quad \text{ for }u,v\in \L_S,$$
where $\bar{u}$ denotes the image in $\L_S/H$ of $u$.

Now, taking $H$ to be subgroups $\ker (\chi)$, $\ker (\tilde{\chi}_B)$, and $\ker (\tilde{\chi}_D)$, we get Lie algebras
\begin{eqnarray}
\L_S^A:=\L_S/\ker (\chi),
\end{eqnarray}
$\L_S^B:=\L_S/\ker (\tilde{\chi}_B)$, and $\L_S^D:=\L_S/\ker(\tilde{\chi}_D)$.
Combining Theorems \ref{cov-A} and \ref{thm-4.9} we immediately obtain:

\bt{A-res}
Let $\ell\in\C$. Then for any restricted $\widehat{A}_{S}$-module $W$ of level $\ell$,
there exists an $(S,\chi)$-equivariant quasi
$V_{\widehat{\L_{S}^A}}(\ell,0)$-module structure $Y_{W}(\cdot, x)$ on $W$,
which is uniquely determined by
$$Y_{W}(\bar{L}_{\alpha,\beta},x)=\chi(\beta)A_{\alpha}(\chi(\beta)x)\quad \mbox{for}\;\alpha,\beta\in S. $$
On the other hand, for any $(S,\chi)$-equivariant quasi $V_{\widehat{\L_{S}^A}}(\ell,0)$-module
$(W, Y_{W})$, $W$ is a restricted $\widehat{A}_{S}$-module of level $\ell$
with
$$A_{\alpha}(x)=Y_{W}(\bar{L}_{\alpha,0},x)\quad \mbox{for}\;\alpha\in S.$$
\et


 With Proposition \ref{cov-B} and Theorem  \ref{thm-4.9} we immediately have:

\bt{BC-res}
Let $\ell\in\C$. Then for any restricted $\widehat{B}_{S}$-module $W$ of level $\ell$,
there exists an $(\tilde{S}_B,\tilde{\chi}_{B})$-equivariant quasi
$V_{\widehat{\L_S^B}}(2\ell,0)$-module structure $Y_{W}(\cdot, x)$ on $W$,
which is uniquely determined by
$$Y_{W}(\bar{L}_{\alpha,\beta},x)=\chi(\beta)B_{\alpha}(\chi(\beta)x)\quad \mbox{for}\;\;\alpha,\beta\in S.$$
On the other hand, for every $(\tilde{S}_B,\tilde{\chi}_{B})$-equivariant
quasi $V_{\widehat{\L_{S}^B}}(2\ell,0)$-module
$(W, Y_{W})$, $W$ is a restricted $\widehat{B}_{S}$-module of level $\ell$
with
$$B_{\alpha}(x)=Y_{W}(\bar{L}_{\alpha,0},x)\quad \mbox{for }\alpha\in S.$$
\et

Combining Proposition \ref{cov-D2} with Theorem  \ref{thm-4.9}  we immediately have:

\bt{D-res}
Let $\ell\in\C$. Then for any restricted $\widehat{D}_{S}$-module $W$ of level $\ell$,
there exists an $(\tilde{S}_D,\tilde{\chi}_D)$-equivariant quasi
$V_{\widehat{\L^{D}_{S}}}(2\ell,0)$-module structure $Y_{W}(\cdot, x)$ on $W$,
which is uniquely determined by
$$Y_{W}(\bar{L}_{\alpha,\beta},x)=\chi(\beta-\alpha)D_{\alpha}(\chi(\beta)x)\quad \mbox{for}\;\alpha,\beta\in S.$$
On the other hand, for any $(\tilde{S}_D,\tilde{\chi}_D)$-equivariant
 quasi $V_{\widehat{\L^{D}_{S}}}(2\ell,0)$-module
$(W, Y_{W})$,  $W$ is a restricted $\widehat{D}_{S}$-module of level $\ell$ with
$$D_{\alpha}(x)=\chi(\alpha)Y_{W}(\bar{L}_{\alpha,0},x)\quad \mbox{for}\;\alpha\in S.$$
\et

\bt{D-res-S}
Let $\ell\in\C$. Then for any restricted $\widehat{D}_{S}$-module $W$ of level $\ell$,
there exists an $(S,\chi)$-equivariant quasi
$V_{\widehat{\L^{\tau}_{S}}}(\ell,0)$-module structure $Y_{W}(\cdot, x)$ on $W$,
which is uniquely determined by
$$Y_{W}(L^{\tau}_{\alpha,\beta},x)=\chi(\beta-\alpha)D_{\alpha}(\chi(\beta)x)\quad \mbox{for}\;\alpha,\beta\in S.$$
On the other hand, for any $(S,\chi)$-equivariant
 quasi $V_{\widehat{\L^{\tau}_{S}}}(\ell,0)$-module
$(W, Y_{W})$,  $W$ is a restricted $\widehat{D}_{S}$-module of level $\ell$ with
$$D_{\alpha}(x)=\chi(\alpha)Y_{W}(L^{\tau}_{\alpha,0},x)\quad \mbox{for}\;\alpha\in S.$$
\et


\section{Lie algebras $\widehat{X}_S$ and affine Kac-Moody algebras}

In this section, we study Lie algebras $\widehat{X}_S$ with $S$ a finite cyclic group
and we show that they are essentially affine Kac-Moody algebras.
More specifically, if $|S|=2l+1$, we show that
$\widehat{A}_S$ is isomorphic to the affine Lie algebra $\widehat{\mathfrak{gl}}_{2l+1}$,
$\widehat{B}_S$ is essentially isomorphic to
the sum of the twisted affine Kac-Moody algebra $\g'(A_{2l}^{(2)})$ and a Heisenberg algebra,
$\widehat{D}_S$ is isomorphic to the affine Kac-Moody algebra $\g'(B_{l}^{(1)})$.
If $|S|=2l$, we show that there exist natural surjective Lie algebra homomorphisms
from $\widehat{X}_S$ to the loop algebra $L(\mathfrak{gl}_{l}(\C))$,
or the loop algebras of classical orthogonal Lie algebras.

\subsection{Lie algebras $\widehat{X}_S$ with $S=\Z_{2l+1}$}

Let $S$ be an abelian group, $\chi:S\rightarrow\C^{\times}$ a linear character, both of which are fixed throughout this section.
First, we give a connection of $\L_{S}$ with a matrix algebra.
Let $\mathfrak{gl}_{S}$ denote the associative algebra with a designated basis
$\{E_{\alpha,\beta}\ |\  \alpha,\beta\in S\}$, where
\begin{eqnarray*}
E_{\alpha,\beta}\cdot E_{\mu,\nu}=\delta_{\beta,\mu}E_{\alpha,\nu}\
\   \  \mbox{ for }\alpha,\beta,\mu,\nu\in S.
\end{eqnarray*}
Note that if $|S|=n$, $\mathfrak{gl}_{S}$ is simply (isomorphic to) the algebra $M(n,\C)$ of
 $n\times n$ complex matrices.
 Equip $\mathfrak{gl}_{S}$ with a bilinear form $\<\cdot,\cdot\>$ defined by
\begin{eqnarray}\label{eq:bil-E}
\<E_{\alpha,\beta}, E_{\mu,\nu}\>=\delta_{\alpha,\nu}\delta_{\beta,\mu}
\;\;\;\;\mbox{for}\;\;\alpha,\beta,\mu,\nu\in S.
\end{eqnarray}
This bilinear form is non-degenerate, symmetric, and associative.
In particular, $\<\cdot,\cdot\>$ on $\mathfrak{gl}_{S}$ viewed as a Lie algebra is invariant.

Note that group $S$ acts on the associative algebra $\mathfrak{gl}_{S}$
as an automorphism group with $\gamma\in S$ acting as $\sigma_{\gamma}$, where
$$\sigma_{\gamma}(E_{\alpha,\beta})=E_{\alpha+\gamma,\beta+\gamma}\ \
\text{for }\alpha,\beta\in S.$$
On the other hand, we have  an order $2$ Lie algebra automorphism $\tau$ of $\mathfrak{gl}_{S}$ defined by
 \begin{eqnarray}
  \tau(E_{\alpha,\beta})=-E_{\beta,\alpha}\quad \text{ for }\alpha,\beta\in S.
  \end{eqnarray}
  Recall that $S^{0}=\{ \alpha\in S\ |\ 2\alpha=0\}$.

The following gives a connection of $\L_{S}$ with the matrix algebra $\mathfrak{gl}_S$:

\bp{LAS}
Let $S$ be an abelian group.
Define a linear map $\pi: \L_{S}\rightarrow \mathfrak{gl}_{S}$ by
\begin{eqnarray}
 \pi(L_{\alpha,\beta})=E_{\alpha+\beta,\beta-\alpha}\quad \text{for }\alpha,\beta\in S.
 \end{eqnarray}
 Then $\pi$ is a homomorphism of associative algebras, which commutes with the actions of $S$ and $\tau$.
Furthermore,  if $S^{0}=0$, then $\pi$ is one-to-one and preserves the bilinear forms.
 If $S^{0}=0$ and $S=2S$, then $\pi$ is an isomorphism.
Especially, if $S$ is finite of an odd order, then $\pi$ is an algebra isomorphism.
\ep

\begin{proof}  For $\alpha,\beta\in S$, set $G_{\alpha,\beta}=E_{\alpha+\beta,\beta-\alpha}$.
We see that for $\alpha,\beta,\mu,\nu\in S$, $G_{\alpha,\beta}=G_{\mu,\nu}$ if and only if
$\alpha-\mu=\beta-\nu\in S^0$ and we have
\begin{eqnarray*}
&G_{\alpha,\beta}\cdot G_{\mu,\nu}=\delta_{\beta-\alpha,\mu+\nu}G_{\alpha+\mu,\alpha+\nu}
=\delta_{\alpha+\mu,\beta-\nu}G_{\alpha+\mu,\alpha+\nu},\\
&\<G_{\alpha,\beta}, G_{\mu,\nu}\>=\delta_{\alpha+\mu,\beta-\nu}\delta_{2(\beta-\nu),0}
\end{eqnarray*}
for $\alpha,\beta,\mu,\nu\in S$. Then the first assertion follows easily.

Assume $S^{0}=0$. Then $\pi(L_{\alpha,\beta})=\pi(L_{\mu,\nu})$ if and only if $(\alpha,\beta)=(\mu,\nu)$.
That is, $\pi$ is one-to-one from the basis $\{ L_{\alpha,\beta}\ | \ \alpha,\beta\in S\}$ of $\L_S$ to
the basis $\{ E_{\alpha,\beta}\ | \ \alpha,\beta\in S\}$ of $\mathfrak{gl}_S$.
It follows that $\pi$ is a one-to-one algebra homomorphism.
With $S^0=0$, we get
\begin{eqnarray*}\label{eq:bil-G}
\<G_{\alpha,\beta}, G_{\mu,\nu}\>=\delta_{\alpha+\mu,0}\delta_{\beta,\nu}
\;\;\;\;\mbox{for}\;\;\alpha,\beta,\mu,\nu\in S,
\end{eqnarray*}
showing that $\pi$ preserves the bilinear forms.

Last, assume that $S$ is finite of an odd order.  In this case, we have $S^0=0$ and $S=2S$.
Then $\pi$ is a bijection between the aforementioned bases of $\L_S$ and $\mathfrak{gl}_S$.
Therefore, $\pi$ is an algebra isomorphism.
\end{proof}

\bc{coro-odd}
Assume that $S$ is finite of odd order.
Then Lie algebra $\widehat{A}_{S}$ is isomorphic to $\widehat{\mathfrak{gl}_{S}}[S]$.
\ec

 \begin{proof} As $S$ is a finite abelian group of odd order, we have $S^{0}=0$ and $2S=S$.
By Proposition \ref{LAS} we get $\L_S\overset{\pi}{\simeq} \mathfrak{gl}_{S}$.
Then using Theorem \ref{cov-A} we obtain
\begin{eqnarray}
\widehat{A}_{S}\overset{\psi_A}{\cong}\widehat{\L_{S}}[S]\overset{\pi}{\cong}\widehat{\mathfrak{gl}_S}[S],
 \end{eqnarray}
as desired.
\end{proof}

Assume $S$ is finite. We have $\mathfrak{gl}_{S}=\mathfrak{sl}_{S}\oplus \C I,$
a direct sum of Lie algebras, where
$\mathfrak{sl}_{S}=\{ X\in \mathfrak{gl}_{S}\ |\ \text{tr}(X)=0\}$
and $I=\sum\limits_{\alpha\in S}E_{\alpha,\alpha}$ is the identity element of $\mathfrak{gl}_{S}$.
Note that the bilinear form, restricted to $\mathfrak{sl}_{S}$, coincides with the normalized killing form
with the squared norm of all (long) roots equal to $2$.

On the other hand, we have $\langle I, I\rangle=|S|$. Set $\mathfrak{h}=\C I$ and form an affine Lie algebra
$\widehat{\mathfrak{h}}=\mathfrak{h}\otimes \C[t,t^{-1}]\oplus\C{\bf k}$. Note that $\tau (I)=-I$.
Define an automorphism $\hat{\tau}$ of $\widehat{\mathfrak{h}}$ by
\begin{eqnarray}\label{def-hat-tau}
\hat{\tau}({\bf k})={\bf k},\ \ \ \hat{\tau}(I\otimes t^m)=-(-1)^m(I\otimes t^m)\quad \text{ for }m\in \Z.
\end{eqnarray}
The $\hat{\tau}$-invariant subalgebra is given explicitly as
\begin{eqnarray}\label{twisted-Heisenberg}
(\widehat{\mathfrak{h}})^{\hat{\tau}}=\mathfrak{h}\otimes t\C[t^2,t^{-2}]+\C{\bf k}.
\end{eqnarray}


Recall from \cite{kac} that for a generalized Cartan matrix $A$,
$\g'(A)$ denotes the derived subalgebra of the Kac-Moody Lie algebra $\g(A)$.
Now, we are ready to present the first main result.

\bt{|S|=odd}
Let $S$ be a finite cyclic group of order $N:=2l+1$ with $l\ge 1$. Then
\begin{itemize}
\item[(1)] Lie algebra $\widehat{A}_{S}$ is isomorphic to the affine Lie algebra $\widehat{\mathfrak{gl}_{N}}$.
\item[(2)] Lie algebra $\widehat{B}_{S}$ is isomorphic to $(\widehat{\mathfrak{gl}_{N}})^{\hat{\tau}}$. Furthermore,
 $(\widehat{\mathfrak{gl}_{N}})^{\hat{\tau}}$ is isomorphic to the quotient algebra of the direct sum
 $\g'(A_{2l}^{(2)})\oplus (\widehat{\mathfrak{h}})^{\hat{\tau}}$ modulo the one-dimensional central ideal $\C ({\bf k},-{\bf k})$.
 \item [(3)]Lie algebra $\widehat{D}_{S}$ is isomorphic to the affine Kac-Moody algebra $\g'(B_{l}^{(1)})$.
\end{itemize}
\et

\begin{proof} Assume $S=\Z_{2l+1}$ with the linear character $\chi$ defined by $\chi(n)=e^{2\pi i n/N}$ for $n\in \Z$.
Set $\sigma=\sigma_{[1]}$, an automorphism of associative algebra $\mathfrak{gl}_{S}$ with
\begin{eqnarray}
\sigma(E_{\alpha,\beta})=E_{\alpha+[1],\beta+[1]}\quad \text{ for }\alpha,\beta\in S, \text{where }[1]\in \Z_{2l+1}.
\end{eqnarray}
Then $S=\langle \sigma\rangle$ is an automorphism group of Lie algebra $\mathfrak{gl}_{S}$.
View $S$ as an automorphism group of the affine Lie algebra $\widehat{\mathfrak{gl}_{S}}$
by using $\chi$ as in Remark \ref{rem-affine-automorphism}.
Using Corollary \ref{coro-odd} and Remark 4.5 of \cite{Li2}, we have
\begin{eqnarray}
\widehat{A}_{S}\overset{\psi_A}{\cong}\widehat{\L_{S}}[S]\overset{\pi}{\cong}\widehat{\mathfrak{gl}_S}[S]
\cong\widehat{\mathfrak{gl}_S}^{S}=\widehat{\mathfrak{gl}_S}^{\sigma}
=\widehat{\mathfrak{gl}_N}^{\sigma}.
 \end{eqnarray}
Set
$$P=
  \begin{pmatrix}
    0 & 0 & 0& \cdots& 0&1\\
    1 & 0 &0&\cdots &0&0\\
    0 & 1& 0 &\cdots & 0&0\\
\vdots&\vdots&\vdots&\vdots&\vdots&\vdots\\
    0 & 0 &0&\cdots &1&0\\
  \end{pmatrix}\in \mathfrak{gl}_N.$$
Then
$\sigma(X)=PXP^{-1}$ for $X\in \mathfrak{gl}_N$.
Let $\log P\in \mathfrak{gl}_N$ such that $e^{\log P}=P$.
 Then
 \begin{eqnarray}
 \sigma(X)=PXP^{-1}=e^{\log P}Xe^{-\log P} =e^{{\rm ad}(\log P)}X
 \end{eqnarray}
for $X\in \mathfrak{gl}_N$ viewed as a Lie algebra. Set
$$h=\log P-\frac{1}{N}{\rm tr}(\log P)I\in  \mathfrak{sl}_N.$$
 We see that $\sigma=e^{{\rm ad}\, h}$ on $\mathfrak{gl}_N$.
 As $\sigma$ is of finite order,  ${\rm ad}\,h$ must be semisimple.
 Then $\widehat{\mathfrak{gl}_N}^{\sigma}\simeq \widehat{\mathfrak{gl}_N}$.
 Therefore, we have
$\widehat{A}_{S}\cong\widehat{\mathfrak{gl}_S}^{\sigma}\cong\widehat{\mathfrak{gl}_N}.$

For Lie algebra $\widehat{B}_{S}$, similarly, using Proposition \ref{cov-B}
we have
\begin{eqnarray}
\widehat{B}_{S}\cong\widehat{\L_{S}}[\tilde{S}_{B}]\cong
 \left(\widehat{\L_{S}}^{\sigma}\right)^{\hat{\tau}}=
 \left(\widehat{\mathfrak{gl}_{S}}^{\sigma}\right)^{\hat{\tau}}\cong
 \left(\widehat{\mathfrak{gl}_{S}}\right)^{\hat{\tau}}.
 \end{eqnarray}
Note that $\mathfrak{gl}_{S}=\mathfrak{sl}_{S}\oplus \mathfrak{h}$ and
$\tau$ preserves both $\mathfrak{sl}_{S}$ and $\mathfrak{h}$. Then
$\hat{\tau}$ is an automorphism of $\widehat{\mathfrak{sl}_{S}}$ and
$\widehat{\mathfrak{h}}$. In particular, this $\hat{\tau}$ on $\widehat{\mathfrak{h}}$ agrees with
the one defined in (\ref{def-hat-tau}).  Explicitly given in (\ref{twisted-Heisenberg}),  $(\widehat{\mathfrak{h}})^{\hat{\tau}}$
is an ideal of $(\widehat{\mathfrak{gl}_{S}})^{\hat{\tau}}$.

Now, consider $(\widehat{\mathfrak{sl}_{S}})^{\hat{\tau}}$, which is also an ideal of $(\widehat{\mathfrak{gl}_{S}})^{\hat{\tau}}$.
 Recall that $\tau(E_{ij})=-E_{ji}$ for $1\le i,j\le N$.
 Let $\mu:\mathfrak{sl}_{N}\rightarrow\mathfrak{sl}_{N}$ be the linear map
 with $\mu(E_{ij})=-E_{N+1-j,N+1-i}$ for $1\le i,j\le N$.
We see that $\mu$ is the standard Dynkin diagram automorphism of $\mathfrak{sl}_{N}$.
Define an automorphism $\alpha$ of the associative algebra $\mathfrak{gl}_{N}$ by
$$\alpha(E_{ij})=E_{N+1-i,N+1-j}\quad \text{ for }1\le i,j\le N.$$
 Set $Q=\sum\limits_{i=1}^{N}E_{i,N+1-i}$. Then $Q^{2}=I$ and $ \alpha(X)=QXQ^{-1}$ for $X\in \mathfrak{gl}_{N}$.
 As before, there exists a matrix $h'\in \mathfrak{gl}_{N}$
such that ${\rm ad}\, h'$ is semisimple and
 \begin{eqnarray*}
 \alpha(X)=QXQ^{-1}=e^{{\rm ad}h'}(X)\quad \text{ for }X\in\mathfrak{gl}_{N}.
 \end{eqnarray*}
 Thus $\tau=\mu\alpha=\mu e^{{\rm ad}h'}$.
  Consequently, $\widehat{\mathfrak{sl}_{S}}^{\hat{\tau}}$ is isomorphic
 to the twisted affine Kac-Moody algebra $\g'(A_{2l}^{(2)})$,
proving the assertion (2).

Lastly, consider Lie algebra $\widehat{D}_{S}$. Combining Propositions \ref{cov-D2} and  \ref{LAS} we get
 \begin{eqnarray*}
 \widehat{D}_{S}\cong\widehat{\L_{S}^{\tau}}[S] \cong\widehat{\mathfrak{gl}_{S}^{\tau}}[S].
 \end{eqnarray*}
Furthermore,  we have
$\widehat{\mathfrak{gl}_{S}^{\tau}}[S]\cong \widehat{\mathfrak{gl}_{S}^{\tau}}^{S}= \widehat{\mathfrak{gl}_{S}^{\tau}}^{\sigma}$.
 Recall that $\sigma=e^{{\rm ad}\, h}$ on $\mathfrak{gl}_S$ and $\sigma\tau=\tau\sigma$.
   It follows that $\widehat{\mathfrak{gl}_{S}^{\tau}}^{\sigma}\cong \widehat{\mathfrak{gl}_{S}^{\tau}}$.
From \cite{gltw2} (Corollary 3.6),
 $\mathfrak{gl}_{S}^{\tau}$ is isomorphic to the classical simple Lie algebra of type $B_{l}$.
 Consequently, $\widehat{D}_{S}$ is isomorphic to the affine Kac-Moody algebra $\g'(B_{l}^{(1)})$.
\end{proof}

\subsection{The Lie algebra $\widehat{A}_{S}$ with $S=\Z_{2l}$}
Here, we study Lie algebra $\widehat{A}_{S}$ for $S=\Z_{2l}$.
As the main result, we prove that there exists a natural surjective homomorphism from $\widehat{A}_{S}$
to the loop algebra $L(\mathfrak{gl}_{l}(\C))$.

Let $S=\Z_{2l}$ with $l$ a positive integer, and fix a faithful linear character
 \begin{eqnarray}
 \chi:\  S=\Z_{2l}\rightarrow \C^{\times}; \quad \chi([m])=e^{m\frac{2\pi i}{2l}}\  \quad  \text{for }m\in \Z.
 \end{eqnarray}
 We also consider $\chi$ as a group homomorphism from $\Z$ to $\C^{\times}$ naturally.
Note that $\chi(2)$ is a primitive $l$-th root of unity and $\chi(l)=-1$.

We shall need a particular basis for $\mathfrak{gl}_{l}(\C)$.
Introduce $l\times l$ matrices
$$P=
  \begin{pmatrix}
    0 & 1 & 0&\cdots&0\\
    0 & 0 & 1&\cdots&0\\
  \vdots&\vdots&\vdots&\vdots&\vdots\\
  0 & 0 & 0&\cdots&1\\
  1 & 0 & 0&\cdots&0\\
  \end{pmatrix},
 \quad \   \
Q=
  \begin{pmatrix}
    \chi(2) & 0 & \cdots&0\\
    0 & \chi(4) & \cdots&0\\
    \vdots & \vdots & \vdots&\vdots\\
    0 & 0 & \cdots&\chi(2l)\\
  \end{pmatrix},
$$
which satisfy the following relations:
\begin{eqnarray}\label{P,Q}
P^{l}=I,\  \ Q^{l}=I,\  \  PQ=\chi(2)QP.
\end{eqnarray}
Furthermore, we have
\begin{eqnarray}\label{P-m-Q-r}
P^{m}Q^{n}=\chi(2mn)Q^{n}P^{m} \;\;\;\;
\mbox{for } m,  n\in\Z.
\end{eqnarray}
It is well known (cf. \cite{Gao}, Corollary 1.12) that $P^{m}Q^{n}$ for $0\le m,n\le l-1$
form a basis for the matrix algebra $M(l,\C)$.
Furthermore, $M(l,\C)$ is isomorphic to the associative algebra generated by two (abstract) elements $P,Q$,
subject to the relation (\ref{P,Q}).
Using this we immediately have:

\bl{lem-anti-theta}
There exists an anti-automorphism $\theta$ of $M(l,\C)$ uniquely determined by
\begin{eqnarray}
\theta(P)=P, \ \ \theta(Q)=Q^{-1}.
\end{eqnarray}
\el

For $r, m\in\Z$, set
\begin{eqnarray}\label{a-r-m}
a_{r,m}=\chi(rm)Q^{r}P^{m}
=\chi(-mr)P^{m}Q^{r}\in M(l,\C).
\end{eqnarray}
The following is an immediate consequence:

\bc{coro-negative-anti}
The linear operator $-\theta$ is a Lie algebra automorphism of $\mathfrak{gl}_{l}(\C)$ with
\begin{eqnarray}
(-\theta)(a_{r,m})=-a_{-r,m}\ \ \text{ for }m,r\in \Z.
\end{eqnarray}
\ec

\bd{def-eta}
{\em Define two algebra automorphisms $\eta_1,\eta_2$ of $M(l,\C)$ by
\begin{eqnarray}
\eta_1(X)=Q^{-1}XQ,\quad \eta_2(X)=PXP^{-1}\quad \text{ for }X\in M(l,\C).
\end{eqnarray}}
\ed

In terms of matrices $a_{r,m}$ with $r,m\in \Z$, we have
\begin{eqnarray}
\eta_1(a_{r,m})=\chi(2m)a_{r,m},\ \ \eta_2(a_{r,m})=\chi(2r)a_{r,m}\quad \text{ for }m,r\in \Z.
\end{eqnarray}
For the general linear Lie algebra $\mathfrak{gl}_{l}(\C)$,
 $\{ a_{r,m}\ | \ 0\le r, m\le l-1\}$ is a basis, where
\begin{eqnarray}\label{eq:A}
 [a_{r,m},a_{s,n}]
=a_{r,m} a_{s,n}-a_{s,n}a_{r,m}
=(\chi(ms-nr)-\chi(nr-ms))a_{r+s,m+n}
\end{eqnarray}
for $r,s, m,n\in\Z$.
The following are obvious relations:
\begin{eqnarray}\label{eq:a-rel}
a_{r,m+l}=(-1)^{r}a_{r,m}, \quad
a_{r+l,m}=(-1)^{m}a_{r,m}
\end{eqnarray}
for $r, m\in\Z$. (Recall that $\chi(l)=-1$.)




\bd{def-bar-eta-1}
{\em Define an automorphism $\bar{\eta}_1$ of the loop algebra $L(\mathfrak{gl}_{l}(\C))$ by
\begin{eqnarray}
 \bar{\eta}_1(a\otimes t^{m})=\chi(2)^{-m}(\eta_1(a)\otimes t^{m})\quad  \mbox{for }a\in \mathfrak{gl}_{l}(\C),\ m\in\Z.
\end{eqnarray}}
\ed

It is straightforward to show that the $\bar{\eta}_1$-invariant
subalgebra $L(\mathfrak{gl}_{l}(\C))^{\bar{\eta}_1}$ has a basis
$$\{a_{r,m}\otimes t^{m}\ | \ 0\leq r\leq l-1,\  m\in\Z\}.$$

For convenience, for $r,m\in\Z$,
we write $A_{[r],m}=A_{r,m}$ for $[r]\in \Z_{2l}$, $r,m\in\Z$.
Now, we are ready to give a natural connection of $\widehat{A}_{S}$ with $L(\mathfrak{gl}_{l}(\C))$.
We have:

\bt{aS}
Assume $S=\Z_{2l}$. Define a linear map $\theta_A:\widehat{A}_{S}\rightarrow
L(\mathfrak{gl}_{l}(\C))^{\bar{\eta}_1}$ by
$$
\theta_A({\bf c})=0,\;\;
\theta_A(A_{r,m})=a_{r,m}\otimes t^{m}\  \  \text{ for } r,  m\in\Z.$$
Then $\theta_A$ is a Lie algebra epimorphism  from $\widehat{A}_{S}$ onto
$L(\mathfrak{gl}_{l}(\C))^{\bar{\eta}_1}$ with
\begin{eqnarray}
\ker \theta_A={\rm span}(\{{\bf c}\}\cup \{A_{r+l,m}-(-1)^{m}A_{r,m} \ | \  r,m\in\Z\}).
\end{eqnarray}
Furthermore, we have
\begin{eqnarray}
L(\mathfrak{gl}_{l}(\C))^{\bar{\eta}_1} \cong L(\mathfrak{gl}_{l}(\C)).
\end{eqnarray}
\et

\begin{proof} First of all, $\theta_A$ is well defined as $a_{r+2l,m}=a_{r,m}$ for $r,m\in \Z$.
Since $L(\mathfrak{gl}_{l}(\C))^{\bar{\eta}_1}$ is spanned by $a_{r,m}\otimes t^{m}$ for $r, m\in\Z$,
$\theta_A$ is surjective.
It follows from (\ref{Abracket}) and (\ref{eq:A}) that $\theta_A$ is a Lie algebra homomorphism.
Let $J_A$ denote the span. As $a_{r+l,m}=(-1)^ma_{r,m}$ for $r,m\in \Z$, we have $J_A\subset \ker (\theta_A)$.
Furthermore, using the aforementioned basis of $L(\mathfrak{gl}_{l}(\C))^{\bar{\eta}_1}$
we conclude $\ker(\theta_A)=J_{A}$.

For the last assertion, set
\begin{eqnarray*}\label{eq:H}
h'= 2\pi i {\rm I}-\frac{2\pi i}{l}{\rm diag}(1,2,\dots,l)\in M(l,\C),
\end{eqnarray*}
where ${\rm diag}(1,2,\dots,l)$ denotes the indicated diagonal matrix.
We have
\begin{eqnarray*}
Q^{-1}=e^{h'}\;\;\;\mbox{and}\;\;\;\eta_1(X)=e^{{\rm ad} h'}(X)
\quad{for}\;X\in \mathfrak{gl}_l(\C).
\end{eqnarray*}
Furthermore, set
\begin{eqnarray}\label{eq:h}
h =h'-\frac{1}{l}{\rm tr}(h')I
=\frac{l+1}{l}\pi i I-\frac{2\pi i}{l}{\rm diag}(1,2,\dots,l)\in \mathfrak{sl}_l(\C).
\end{eqnarray}
We see that $\eta_1(X)=e^{{\rm ad} h}(X)$ for $X\in \mathfrak{gl}_{l}(\C)$. Then we have
$L(\mathfrak{gl}_{l}(\C))^{\bar{\eta}_1} \cong L(\mathfrak{gl}_{l}(\C)).$
(This also follows from Theorem 9.1 of \cite{ABP} as $\eta_1$ clearly fixes point-wise the standard
Cartan subalgebra of $\mathfrak{gl}_{l}(\C)$.)
\end{proof}

\subsection{The Lie algebra $\widehat{C}_{S}$ with $S=\Z_{2l}$}

Note that for $r, m\in \Z$, we have
$$\theta_A(C_{r,m})=(a_{r,m}-(-1)^m\chi(2r)a_{-r,m})\otimes t^m\in L(\mathfrak{gl}_{l}(\C))^{\bar{\eta}_1}.$$
Motivated by this, for $r, m\in \Z$, set
\begin{eqnarray}\label{def-c-r-m}
c_{r,m}=a_{r,m}-(-1)^{m}\chi(2r)a_{-r,m}\in \mathfrak{gl}_{l}(\C).
\end{eqnarray}
Then $\theta_A(\widehat{C}_S)$ is linearly spanned by $c_{r,m}\otimes t^m$ for $r,m\in \Z$.
It follows that $c_{r,m}$ for $r,m\in \Z$ linearly span a Lie subalgebra, denoted by $c_S$.

From (\ref{eq:a-rel}) we also have
\begin{eqnarray}\label{eq:ca-rel}
c_{r,m+l}=(-1)^{r}(a_{r,m}-(-1)^{m+l}\chi(2r)a_{-r,m})
\end{eqnarray}
for $r,m\in \Z$. In case that $l$ is odd, combining
(\ref{def-c-r-m}) and (\ref{eq:ca-rel}) we obtain
$$a_{r,m}=\frac{1}{2}(c_{r,m}+(-1)^{r}c_{r,m+l})\quad \text{ for  }r, m\in\Z.$$
Thus we have:

\bl{l=odd}
Assume that $l$ is a positive odd integer.
Then $\mathfrak{gl}_{l}(\C)=c_S$, i.e., $\mathfrak{gl}_{l}(\C)$ is linearly spanned by $c_{r,m}$ for $r,m\in \Z$.
\el

\bl{rel-c}
The following relations hold for $r, m\in\Z$:
 \begin{eqnarray}\label{three-conditions-c}
 c_{rl,2m}=0,\quad
 c_{r+l,m}=(-1)^{m}c_{r,m},\quad c_{r,m+2l}=c_{r,m}.
 \end{eqnarray}
 If $l$ is a positive even integer, the following relations also hold for $r, m\in\Z$:
\begin{eqnarray}\label{c-relations}
c_{r,m+l}=(-1)^{r}c_{r,m}, \quad  c_{r+\frac{l}{2},m}=\chi(2r)c_{\frac{l}{2}-r,m}.
\end{eqnarray}
\el

\begin{proof} The relations in (\ref{three-conditions-c}) follow straightforwardly from (\ref{def-c-r-m}) and
 (\ref{eq:a-rel}).
 Assume that $l$ is even. For $r,m\in\Z$, we have
 \begin{eqnarray*}
 \begin{aligned}
c_{r,m+l}
&=(-1)^r\left(a_{r,m}-(-1)^m\chi(2r)a_{-r,m}\right)=(-1)^{r}c_{r,m},\\
c_{r+\frac{l}{2},m}
&=a_{r+\frac{l}{2},m}-(-1)^{m}\chi(2r+l)a_{-r-\frac{l}{2},m} \\
&=(-1)^ma_{r-\frac{l}{2},m}+\chi(2r)a_{-r+\frac{l}{2},m}\\
&=\chi(2r)c_{-r+\frac{l}{2},m},
\end{aligned}
\end{eqnarray*}
recalling that $\chi(l)=-1$.
\end{proof}

Next, we consider the case with $l$ a positive even integer.

\bl{def-tau-c}
Assume that $l$ is a positive even integer.
There exists an order-$2$ automorphism $\tau_c$ of Lie algebra
 $\mathfrak{gl}_{l}(\C)$, which is uniquely determined by
\begin{eqnarray}\label{tau-c-def}
\tau_{c}(a_{r,m})=-(-1)^{m}\chi(2r)a_{-r,m}\quad \text{ for }r, m\in\Z.
\end{eqnarray}
Furthermore, we have
\begin{eqnarray}
\theta_A \tau_C=(\tau_c\otimes 1)\theta_A,\ \
\tau_c\eta_1=\eta_1\tau_c,
\end{eqnarray}
where $\tau_c\otimes 1$ is the indicated automorphism of $L(\mathfrak{gl}_{l}(\C))=\mathfrak{gl}_{l}(\C)\otimes \C[t,t^{-1}]$.
\el

\begin{proof} Recall from Definition \ref{def-eta} the algebra automorphisms $\eta_1,\eta_2$ of $M(l,\C)$ with
$$\eta_1(a_{r,m})=\chi(2m)a_{r,m},\quad \eta_2(a_{r,m})=\chi(2r)a_{r,m}\ \text{ for }m,r\in \Z.$$
As $\chi(l)=-1$, we have
$$\eta_1^{l/2}\eta_2(a_{r,m})=\chi(lm)\chi(2r)a_{r,m}=(-1)^m\chi(2r)a_{r,m}.$$
On the other hand, by Corollary \ref{coro-negative-anti}, $-\theta$ is a Lie algebra automorphism  of
$\mathfrak{gl}_{l}(\C)$ with $(-\theta)(a_{r,m})=-a_{-r,m}$ for $r,m\in \Z$.
Then $\tau_c:=(-\theta)\eta_1^{l/2}\eta_2$ is a Lie algebra automorphism of $\mathfrak{gl}_{l}(\C)$
satisfying the required condition. The proof for the ``furthermore'' assertion is straightforward.
\end{proof}

Furthermore, we have the following result:

\bl{Lie-cS-even}
Assume $l$ is a positive even integer. Then
$\mathfrak{gl}_{l}(\C)^{\tau_c}$ is isomorphic to the
classical simple Lie algebra $\mathfrak{sp}(l,\C)$ of type $C_{\frac{l}{2}}$.
\el

\begin{proof} From Lemma \ref{rel-c}, $\mathfrak{gl}_{l}(\C)^{\tau_c}$ is linearly spanned by $c_{0,2j+1}$
for $0\leq j\leq \frac{l}{2}-1$ and $c_{r,m}$ for $1\leq r\leq \frac{l}{2},  0\leq m\leq l-1$.
It is straightforward to show that this spanning set is linearly independent, so that
$\dim \mathfrak{gl}_{l}(\C)^{\tau_c}=\frac{1}{2}l^{2}+\frac{1}{2}l$.
Set
$$T_{1}=
  \begin{pmatrix}
     0 & 0 &\cdots&0&-1\\
    0 & 0 &\cdots&1&0\\
  \vdots&\vdots&\vdots&\vdots&\vdots\\
     0 & -1 &\cdots&0&0\\
  1 & 0 & \cdots&0&0\\
  \end{pmatrix}\in M(l,\C).
$$
It is straightforward to show that
$$T_1^{-1}=-T_1, \ \ T_1PT_1^{-1}=-P^{t},\ \ T_1QT_1^{-1}=\chi(2)Q^{-1},\ \ QP^{t}=\chi(2)P^{t}Q,$$
where $M^t$ stands for the transpose of a square matrix $M$. Using these facts we get
\begin{eqnarray*}
&&T_1a_{r,m}T_1^{-1}=\chi(rm)\chi(2r)(-1)^mQ^{-r}(P^t)^m=\chi(-rm)\chi(2r)(-1)^m(P^t)^mQ^{-r}\\
&&=(-1)^m\chi(2r)a_{-r,m}^t.
\end{eqnarray*}
Furthermore, we have $T_1c_{r,m}T_1^{-1}=-c_{r,m}^t$, or equivalently,
$$ c_{r,m}^{t}T_{1}=-T_{1}c_{r,m}\quad \text{for }r,m\in\Z.$$
As $\dim \mathfrak{gl}_{l}(\C)^{\tau_c}= \frac{l^{2}}{2}+\frac{l}{2}=\dim \mathfrak{sp}(l,\C)$,
we conclude $\mathfrak{gl}_{l}(\C)^{\tau_c}\cong\mathfrak{sp}(l,\C)$.
\end{proof}

Notice that for the automorphism $\eta_1$ of $\mathfrak{gl}_{l}(\C)$, we have
\begin{eqnarray*}
\eta_1(c_{r,m})=\chi(2m)c_{r,m}\quad \text{ for }r, m\in\Z.
\end{eqnarray*}
Then $\eta_1$ is also an automorphism of the Lie algebra $c_{S}$.
It is clear that
$$L(c_{S})^{\bar{\eta}_1}=\span\{c_{r,m}\otimes t^{m}\ | \ 0\leq r\leq l-1,  m\in\Z\}.$$

\bl{c-loop}
If $l$ is an odd positive integer, then $L(c_{S})^{\bar{\eta}_1}\cong L(\mathfrak{gl}_{l}(\C))$.
If $l$ is even, then $L(c_{S})^{\bar{\eta}_1}\cong L(\mathfrak{sp}(l,\C))$.
\el

\begin{proof} Assume $l$ is odd. By Lemma \ref{l=odd}, we have $c_{S}=\mathfrak{gl}_{l}(\C)$.
Then using Theorem \ref{aS} we get
$$L(c_{S})^{\bar{\eta}_1}=L(\mathfrak{gl}_{l}(\C))^{\bar{\eta}_1}\cong L(\mathfrak{gl}_{l}(\C)).$$
Now, assume $l$ is even. By Lemmas \ref{def-tau-c} and \ref{Lie-cS-even}, we have
$c_{S}=\mathfrak{gl}_{l}(\C)^{\tau_c}\cong \mathfrak{sp}(l,\C)$.
Recall the matrix $h$ from (\ref{eq:h}). It is easy to show that $h^{t}T_{1}=-T_{1}h$,
so that $h\in\mathfrak{gl}_{l}(\C)^{\tau_c}$.
As $\eta_1=e^{{\rm ad} h}$ on $\mathfrak{gl}_{l}(\C)^{\tau_c}$,
we obtain (cf. Theorem 9.1 of \cite{ABP})
$$L(c_{S})^{\bar{\eta}_1}=L(\mathfrak{gl}_{l}(\C)^{\tau_c})^{\bar{\eta}_1}\cong L(\mathfrak{gl}_{l}(\C)^{\tau_c})
\cong L(\mathfrak{sp}(l,\C)),$$
as desired.
\end{proof}

For convenience, we formulate the following straightforward result:

\bl{lem-basic}
Let $\sigma$, $\tau$ be order-$2$ linear automorphisms of vector spaces $U$, $V$, respectively,
and let $\psi: U\rightarrow V$ be a surjective linear map such that $\psi\sigma=\tau\psi$. Then
$\psi(U^{\sigma})=V^{\tau}$, $\sigma (\ker \psi)=\ker \psi$, and
\begin{eqnarray}
\ker (\psi|_{U^{\sigma}})=(\ker \psi)\cap U^{\sigma}={\rm span}\{ u+\sigma (u)\ |\ u\in \ker \psi\}.
\end{eqnarray}
\el

With Theorem \ref{aS} and the relation $\theta_A \tau_C=(\tau_c\otimes 1)\theta_A$ (Lemma \ref{def-tau-c}),
using Lemmas \ref{lem-basic}, \ref{Lie-cS-even}, and \ref{c-loop},   we obtain the following result:

\bt{bcS}
Let $S=\Z_{2l}$ with $l$ a positive integer. Then there is a Lie algebra epimorphism $\theta_C$
from $\widehat{C}_{S}$ to the loop algebra $L(c_{S})^{\bar{\eta}_1}$ such that
$$\theta_C({\bf c})=0,\;\;
\theta_C(C_{r,m})=c_{r,m}\otimes t^{m}\quad \text{for } r, m\in\Z$$
and such that
\begin{eqnarray}
\ker \theta_C={\rm span} (\{{\bf c}\}\cup \{C_{r+l,m}-(-1)^{m}C_{r,m}
\ | \ r, m\in\Z \}).
\end{eqnarray}
Furthermore, we have
$$L(c_{S})^{\bar{\eta}_1}\cong L(\mathfrak{gl}_{l}(\C))\ \  \text{ if $l$ is odd},$$
$$L(c_{S})^{\bar{\eta}_1}\cong L(\mathfrak{sp}(l,\C))\ \  \text{ if $l$ is even}.$$
\et

\subsection{The Lie algebra $\widehat{D}_{S}$ with $S=\Z_{2l}$}







Here, we show that there exists a natural Lie algebra epimorphism from $\widehat{D}_{S}$ onto
the loop algebras of classical orthogonal Lie algebras (of type $B$ or type $D$).

For $l=1$, i.e., $S=\Z_{2}$, we have $2\alpha=0$, i.e., $-\alpha=\alpha$ for $\alpha\in S$.
 For $\alpha\in S,m\in \Z$, since $D_{-\alpha,m}=-\chi(-2\alpha)D_{\alpha,m}$ we conclude $D_{\alpha,m}=0$.
Thus  $\widehat{D}_{S}=0$.
In the rest of this subsection, we assume $l>1$.

Recall the algebra automorphism $\eta_2$ of $M(l,\C)$ from Definition \ref{def-eta},
 the Lie algebra automorphism $-\theta$ of $\mathfrak{gl}_{l}(\C)$ from Corollary \ref{coro-negative-anti}, and
the Lie algebra epimorphism $\theta_A: \widehat{A}_S\rightarrow L(\mathfrak{gl}_{l}(\C))^{\bar{\eta}_1}$ from
Theorem \ref{aS}.
 The following is immediate:

\bl{def-tau-d}
Set $\tau_d=(-\theta)\eta_2$. Then $\tau_d$ is an order-$2$ Lie algebra automorphism of $\mathfrak{gl}_{l}(\C)$,
 which is uniquely determined by
\begin{eqnarray}\label{tau-d-def}
\tau_{d}(a_{r,m})=-\chi(2r)a_{-r,m}\quad \text{ for }r, m\in\Z.
\end{eqnarray}
Furthermore, we have
\begin{eqnarray}
\theta_A \tau_D=(\tau_d\otimes 1)\theta_A,\ \
\tau_d\eta_1=\eta_1\tau_d.
\end{eqnarray}
\el


Next, we identify $\mathfrak{gl}_{l}(\C)^{\tau_d}$ (cf. \cite{G-K-L1}, \cite{G-K-L2}).
For $r, m\in \Z$, set
\begin{eqnarray}
d_{r,m}=a_{r,m}+\tau_d(a_{r,m})=a_{r,m}-\chi(2r)a_{-r,m}.
\end{eqnarray}
Then
\begin{eqnarray}
\mathfrak{gl}_{l}(\C)^{\tau_d}={\rm span}\{ d_{r,m}\ |\ r,m\in \Z\}.
 \end{eqnarray}

\bl{d-rel}
The following relations hold for $r,m\in\Z$:
 \begin{eqnarray}\label{eq:c1}
 d_{0,m}=0,\;\;\; d_{r,m+l}=(-1)^{r}d_{r,m},\;\;\;
 d_{r+l,m}=(-1)^{m}d_{r,m}.
 \end{eqnarray}
Furthermore, if $l$ is even, the following relation also hold
for $r,m\in\Z$:
\begin{eqnarray}\label{eq:c2}
d_{r+\frac{l}{2},m}=(-1)^{m}\chi(2r)d_{-r+\frac{l}{2},m}.
\end{eqnarray}
In particular, $d_{\frac{l}{2},m}=0$ for $m$ odd.
 If $l$ is odd, the following relation hold:
\begin{eqnarray}\label{eq:c3}
d_{r+\frac{l+1}{2},m}=(-1)^{m}\chi(2r+1)d_{-r+\frac{l-1}{2},m}
\end{eqnarray}
 for $r,m\in\Z$.
\el

\begin{proof}
The relations in (\ref{eq:c1}) follow from the definition of $d_{r,m}$ and (\ref{eq:a-rel}).
Assume that $l$ is even.  For $r,m\in\Z$, we have
\begin{eqnarray*}
d_{r+\frac{l}{2},m}
&=&a_{r+\frac{l}{2},m}-\chi(2r+l)a_{-r-\frac{l}{2},m}\\
&=&(-1)^{m}(a_{r-\frac{l}{2},m}+\chi(2r)a_{-r+\frac{l}{2},m})\\
&=&(-1)^{m}\chi(2r)(a_{-r+\frac{l}{2},m}-\chi(-2r+l)a_{r-\frac{l}{2},m})\\
&=&(-1)^{m}\chi(2r)d_{-r+\frac{l}{2},m},
\end{eqnarray*}
recalling that $\chi(l)=-1$. In particular, we have $d_{\frac{l}{2},m}=0$ for $m$ odd.

If $l$ is odd, for $r,m\in\Z$ we have
\begin{eqnarray*}
d_{r+\frac{l+1}{2},m}
&=& a_{r+\frac{l+1}{2},m}+\chi(2r+1)a_{-r-\frac{l+1}{2},m}\\
&=& (-1)^{m}(a_{r-\frac{l-1}{2},m}+\chi(2r+1)a_{-r+\frac{l-1}{2},m}) \\
&=&(-1)^{m}\chi(2r+1)(a_{-r+\frac{l-1}{2},m}-\chi(-2r+l-1)a_{-(-r+\frac{l-1}{2}),m})\\
&=&(-1)^{m}\chi(2r+1)d_{-r+\frac{l-1}{2},m},
\end{eqnarray*}
as desired.
\end{proof}

\bl{Lie-dS-odd}
Let $l$ be a positive integer with $l\ge 2$. Then
 $\mathfrak{gl}_{l}(\C)^{\tau_d}$ is isomorphic to
the classical simple Lie algebra $\mathfrak{o}(l,\C)$, which is of type $B_{\frac{l-1}{2}}$ if $l$ is odd,
and of type $D_{\frac{l}{2}}$ if $l$ is even.
\el

\begin{proof} Set
$$T=
  \begin{pmatrix}
    0 & 0 &\cdots&0&1\\
    0 & 0 &\cdots&1&0\\
  \vdots&\vdots&\vdots&\vdots&\vdots\\
  0 & 1 & \cdots&0&0\\
  1 & 0 & \cdots&0&0\\
  \end{pmatrix}\in M(l,\C).
$$
It is straightforward to show that
$$T^{-1}=T,\ \ TQT^{-1}=\chi(2)Q^{-1},\ \ TPT^{-1}=P^t.$$
Using these properties we get
$$ d_{r,m}^{t}T=-Td_{r,m}\quad \text{for }\;r,m\in\Z.$$
Assume $l$ is even. From Lemma \ref{d-rel}, $\mathfrak{gl}_{l}(\C)^{\tau_d}$ is spanned by
$d_{\frac{l}{2},m}$ for $0\leq m\leq l-1$ even,
and $d_{r,m}$ for $1\leq r\leq \frac{l}{2}-1$, $0\leq m\leq l-1$.
It is straightforward to show that the elements are also linearly independent.
 As $\dim \mathfrak{gl}_{l}(\C)^{\tau_d}= \frac{l^{2}}{2}-\frac{l}{2}=\dim \mathfrak{o}(l,\C)$,
it follows that $\mathfrak{gl}_{l}(\C)^{\tau_d}\cong\mathfrak{o}(l,\C)$.
Similarly,
if $l$ is odd, using Lemma \ref{d-rel} we can show that
$d_{r,m}$ for $1\leq r\leq \frac{l-1}{2}$, $0\leq m\leq l-1$ form a basis for $\mathfrak{gl}_{l}(\C)^{\tau_d}$.
Then by comparing the dimensions we conclude that
 $\mathfrak{gl}_{l}(\C)^{\tau_d}\cong \mathfrak{o}(l,\C)$.
\end{proof}

The automorphism $\eta_1$ of $\mathfrak{gl}_{l}(\C)$ gives an automorphism of
 $\mathfrak{gl}_{l}(\C)^{\tau_d}$ with
  $$\eta_1(d_{r,m})=\chi(2m)d_{r,m}\quad\mbox{for}\;r,m\in\Z.$$
Furthermore, $\bar{\eta}_1$ is an automorphism  of the loop algebra $L(\mathfrak{gl}_{l}(\C)^{\tau_d})$ and
$$L(\mathfrak{gl}_{l}(\C)^{\tau_d})^{\bar{\eta}_1}=\span\{d_{r,m}\otimes t^{m}\ | \ 1\leq r\leq l-1, m\in\Z\}.$$
 We see that the matrix $h$ in (\ref{eq:h}) satisfies $h^{t}T=-Th$, so $h\in \mathfrak{gl}_{l}(\C)^{\tau_d}$
and $\eta_1=e^{{\rm ad} h}$ on $\mathfrak{gl}_{l}(\C)^{\tau_d}$.
Just as before, using Lemma \ref{Lie-dS-odd}, we get
\begin{eqnarray}
L(\mathfrak{gl}_{l}(\C)^{\tau_d})^{\bar{\eta}_1}\cong L(\mathfrak{gl}_{l}(\C)^{\tau_d})\cong L(\mathfrak{o}(l,\C)).
\end{eqnarray}


To summarize we have:

\bt{dS}
Assume $S=\Z_{2l}$ with $l\ge 2$. Then there exists a Lie algebra epimorphism $\theta_D$
from $\widehat{D}_{S}$ onto $L(\mathfrak{gl}_{l}(\C)^{\tau_d})^{\bar{\eta}_1}$ such that
$$\theta_D({\bf c})=0,\;\;
\theta_D(D_{r,m})=d_{r,m}\otimes t^{m} \quad \text{for } r,  m\in\Z,$$
where
$$\ker \theta_D={\rm span}(\{ {\bf c}\}\cup \{ D_{r+l,m}-(-1)^mD_{r,m}\ |\ r,m\in \Z\}).$$
Furthermore, we have
\begin{eqnarray*}
L(\mathfrak{gl}_{l}(\C)^{\tau_d})^{\bar{\eta}_1}\cong L(\mathfrak{o}(l,\C)).
\end{eqnarray*}
\et


\end{document}